\documentclass[a4paper,11pt]{article}
\usepackage{amsmath,amssymb,amsthm}
\usepackage{hyperref,url,enumitem,bm,esint}
\usepackage{color}
\usepackage[margin=30mm]{geometry}

\definecolor{rosso}{rgb}{0.8,0,0}
\definecolor{gray}{gray}{0.7}
\def\andreaold #1{{\color{rosso}#1}}

\def\andrea #1{{\color{black}#1}}

\def\andreaold #1{#1}
\newcommand{\x}{\bm{x}}
\newcommand{\A}{\mathcal{A}}
\newcommand{\M}{\mathcal{M}}
\newcommand{\N}{\mathbb{N}}
\newcommand{\0}{\bm{0}}

\newcommand{\R}{\mathbb{R}}
\newcommand{\bet}{\bm{\eta}}

\newcommand{\bY}{\bm{Y}}
\newcommand{\dA}{\, dA}
\newcommand{\bsi}{{\mathcal S}}

\renewcommand{\u}{\bm{u}}
\newcommand{\E}{\mathcal{E}}
\newcommand{\bE}{\bar{\mathcal{E}}}
\newcommand{\C}{\mathcal{C}}
\renewcommand{\div}{\mathrm{div}}
\newcommand{\y}{\bm{y}}
\newcommand{\dx}{\, d \x}
\newcommand{\dt}{\, dt}

\newcommand{\pd}{\partial}
\newcommand{\nx}{\nabla}
\newcommand{\Lx}{\Delta}
\newcommand{\pdnu}{\pd_{\bf n}}
\newcommand{\abs}[1]{\left| #1 \right|}
\newcommand{\no}[1]{\| #1 \|}
\newcommand{\inn}[2]{\langle #1 , #2 \rangle}
\newcommand{\eps}{\varepsilon}
\newcommand{\mean}[1]{\langle #1 \rangle}

\newcommand{\GD}{\Gamma_D}
\newcommand{\GN}{\Gamma_N}
\newcommand{\SD}{\Sigma_D}
\newcommand{\SN}{\Sigma_N}
\newcommand{\g}{\bm g}
\newcommand{\X}{X (\Omega)}
\newcommand{\ph}{\varphi}
\newcommand{\s}{\sigma}

\def\genspazio #1#2#3#4#5{#1^{#2}(#5,#4;#3)}
\def\spazio #1#2#3{\genspazio {#1}{#2}{#3}T0}

\def\L {\spazio L}
\def\H {\spazio H}
\def\Hx #1{H^{#1}(\Omega)}
\def\LOx #1{L^{#1}(\Omega)}

\def\norma #1{\mathopen \| #1\mathclose \|}
\def\d{\delta}

\theoremstyle{plain}
\newtheorem{thm}{Theorem}[]

\newtheorem{remark}{Remark}[section]

\newtheorem{defn}{Definition}[section]
\numberwithin{equation}{section}

\title{On a phase field model of Cahn--Hilliard type for tumour growth with mechanical effects}

\author{Harald Garcke \footnotemark[1] \and Kei Fong Lam \footnotemark[2]
\and Andrea Signori \footnotemark[3]}
\date{\today}

\begin{document}

\maketitle

\renewcommand{\thefootnote}{\fnsymbol{footnote}}
\footnotetext[1]{Fakult{\"a}t f\"ur Mathematik, Universit{\"a}t Regensburg, 93040 Regensburg, Germany
({\tt Harald.Garcke@mathematik.uni-regensburg.de}).}
\footnotetext[2]{Department of Mathematics, The Chinese University of Hong Kong, Shatin, N.T., Hong Kong ({\tt kflam@math.cuhk.edu.hk}).}
\footnotetext[3]
{Dipartimento di Matematica e Applicazioni, Universit\`a di Milano--Bicocca
via Cozzi 55, 20125 Milano, Italy ({\tt andrea.signori02@universitadipavia.it}).}
\begin{abstract}
\noindent
Mechanical effects  have mostly  been   neglected so far  in  phase field tumour  models that are based
on a Cahn--Hilliard approach. In this paper we study a macroscopic mechanical
model for tumour growth in which cell-cell adhesion effects are taken into
account with the help of a Ginzburg--Landau type energy. In the overall model
an equation of
Cahn--Hilliard type is coupled to the system of linear elasticity and
a reaction-diffusion equation for a nutrient concentration.
The highly non-linear coupling between a fourth-order Cahn--Hilliard equation and  the quasi-static
elasticity system lead to new challenges which cannot be dealt within a
gradient flow  setting which was the method of choice for other elastic
Cahn--Hilliard systems. We show existence, uniqueness and regularity results.
In addition, several continuous dependence results with respect to different
topologies are shown. Some of these results give
uniqueness for weak solutions and other results will be helpful for optimal
control problems.   
%\andrea{
%A phase field model of Cahn-Hilliard type  
%describing tumour growth including elastic effects
%and chemotaxis is discussed. 
%Moreover, a mixed Dirichlet-Neumann boundary condition for the 
%stress is postulated to describe the possible presence of some bones which
%would prevent the displacement to vary in that region.
%Well-posedness results for the tumour model is established and 
%continuous dependence results are provided under suitable assumptions
%on the elastic energy (namely, Vegard's law and homogeneous elasticity). 
%Lastly, some asymptotic techniques allow us
%to discuss a variant in case of quasi-static nutrient.}
\end{abstract}


\section{Introduction}
Modelling of tumour growth is one of the challenging frontiers of applied
mathematics. In the last years phase field models for tumour growth have been
studied intensively. Alike classical free boundary models they use a continuum
approach to describe the growth of tumours. However, an advantage to free
boundary models is that phase field models allow for topology changes like break up and
coalescence. In phase field models an order parameter is introduced to
describe the tumour fraction locally in space. In this paper the 
order parameter is denoted by $\varphi$ and it will take the value $+1$ in
regions occupied solely by tumour cells and $-1$ in regions occupied solely 
by healthy
cells. 
As summarised in Lima et al.~\cite{Lima0,Lima} stress effects resulting from tumour
growth severely affect the growth itself. Experimental studies, see
\cite{CTJM, HNLMJ, SMCetal},  show that stresses can inhibit tumour
growth. In this paper we will consider the effect of stresses on the mobility
and on the proliferation rate. In particular, the mobility and the
proliferation rate  will decrease with
increasing stresses.
Often mechanically-coupled models for tumour growth use reaction-diffusion
systems to model proliferation and nutrient diffusion with specific body force
fields that take elastic effects into account, see \cite{Hogea,
  JonesByrneetal, Weisetal}. However, in the recent work of Lima et al.~\cite{Lima0} on
selection, calibration and validation of 
different models of tumour growth, it turned out that phase field methods
taking elastic effects into account are the best modelling approach and in particular
superior to reaction-diffusion models, see the conclusion section of
\cite{Lima0}.
It is the goal of this paper to generalise the model studied in
\cite{Lima0,Lima} and in particular also take nutrient diffusion into account as the latter
definitely will have an important effect on tumour growth.

We will consider  balance equations for the tumour and nutrient concentrations
which will be of parabolic type.
As  diffusion and growth take place on a timescale much larger than that
associated with inertia, we disregard inertial terms and consider instead a
quasi-static approximation. 
We hence consider the following  extension of the phase field tumour model proposed by Lima et al.~\cite{Lima0, Lima}.  
Let $\Omega \subset \R^d$, $d = 2,3$, denote a $C^{1,1}$ or convex bounded domain,
with boundary $\Gamma := \pd \Omega$ and let $\Gamma_D,\Gamma_N\subset \Gamma$ 
are relatively open such that
\begin{align*}
	\Gamma = \overline{\Gamma}_D \cup \overline{\Gamma}_N, 
	\quad 
	\GD \cap \GN = \emptyset, 
	\quad 
	|\GD| > 0,
\end{align*}
where $|\Gamma_D |$ stands for the $(d-1)$-dimensional Hausdorff measure of $\Gamma_D$.  For a fixed but arbitrary time $T > 0$, we consider the system 
\begin{subequations}\label{Lima}
\begin{alignat}{3}
\varphi_t & = \div (m(\varphi, \sigma, \u, \E(\u)) \nx \mu) + U(\varphi, \sigma, \E(\u)) && \text{ in } Q := \Omega \times (0,T), \label{phi} \\
\mu & = - \eps \Lx \varphi + \eps^{-1} \psi'(\varphi)  - \chi \sigma  + W_{,\varphi}(\varphi, \E(\u)) && \text{ in } Q, \label{mu} \\
W_{,\varphi} & = \tfrac{1}{2} (\E(\u) - \bE(\varphi)) :\C'(\varphi)(\E(\u) - \bE(\varphi)) \label{W,phi} \\
\notag & \quad - \C(\varphi) (\E(\u) - \bE(\varphi)) : \bE'(\varphi), \\
\beta \sigma_t & = \Lx \sigma + S(\varphi,\sigma) && \text{ in } Q, \label{sigma} \\
\0 & = \div (W_{,\E}(\varphi, \E(\u))) && \text{ in }Q, \label{u} \\
W_{,\E} & = \C(\varphi)(\E(\u) - \bE(\varphi)), \label{W,E} \\
\varphi(.,0) & = \varphi_0(.), \quad \sigma(.,0) = \sigma_0(.) && \text{ in } \Omega, \label{ic} \\
0 & = \pdnu \varphi = \pdnu \mu, \quad \pdnu \sigma + \kappa (\sigma - \sigma_B) = 0&& \text{ on } \Sigma := \Gamma \times (0,T), \\
\u & = \0  && \text{ on } \SD:= \GD \times (0,T), \label{SD}\\
\label{g} W_{,\E} \, \bold n &= \g
&& \text{ on } \SN:= \GN \times (0,T).
\end{alignat}
\end{subequations}

In the above equations, the primary variables of the model are $\varphi$ (the difference in volume fractions between the tumour and healthy cells), $\mu$ (the associated chemical potential), $\sigma$ (the nutrient concentration) and $\u$ (the displacement). 
Moreover, $\bold n$ indicates the outward unit normal of $\Gamma$ and $\partial_{\bf n}$
stands for the outward normal derivative. The quantity $\E(\u) = \frac{1}{2}(\nx \u + (\nx \u)^t)$ is the symmetric strain tensor, $\psi'$ is the derivative of a double-well function $\psi$ (with the classical example being $\psi(s) = (s^2-1)^2$), $W_{,\varphi}, W_{,\E}$ denote the partial derivatives of the elastic energy $W(\varphi, \E(\u))$ with respect to its arguments, and $\eps > 0$ is a positive parameter associated to the thickness of the interfacial layer.  It is worth noting that, in the phenomena we are interested in, the strain is usually small so that the linearised strain
tensor is used.

Equations \eqref{phi}-\eqref{mu} comprises of a Cahn--Hilliard system
for $(\varphi, \mu)$ with a positive mobility
$m(\varphi, \sigma, \u, \E(\u))$ and a source term
$U(\varphi, \sigma, \E(\u))$ modelling the growth and death of cells.
As one biologically relevant example for the source term $U$ we suggest
\begin{align}\label{defn:U}
U(\varphi, \sigma, \E(\u))) =  \frac{\lambda_p f(\varphi) \sigma}{1+ \abs{W_{,\E}(\varphi, \E(\u))}} - \lambda_a k(\varphi),
\end{align}
for some bounded functions $f(\varphi)$ and $k(\varphi)$. The coefficients $\lambda_p \geq 0$ and $\lambda_a \geq 0$ have the meaning of proliferation and apoptosis rates, respectively.  As discussed in \cite[p.~353]{Byrne}, one should also account for the effects of mechanical interactions in cell growth, such as the stress exerted on the replicating cells by the surrounding environments which leads to a strong dependence of cellular proliferation on mechanical stresses.  The above choice \eqref{defn:U} ensures that as the magnitude of the stress $
\bsi:=W_{,\E}(\varphi, \E(\u))$ increases, the effects of
proliferation are reduced. Other possible forms of $U$  include the
von Mises stress as a stress measure, see \cite{Lima0}, and could also be
used in the theory stated later.

In \eqref{mu}, directed movement of the tumour cells by chemotaxis is captured in the term $-\chi \sigma$, so that $\chi \geq 0$ is a chemotactic sensitivity, and the effects of elastic deformation on the movement of the tumour cells is prescribed by the term $W_{,\varphi}(\varphi, \E(\u))$, whose full expression is given in \eqref{W,phi}.

In \eqref{W,phi}, $\C(\varphi)$ is a symmetric and positive definite
elasticity tensor depending on $\varphi$ and $\bE(\varphi)$ is the
stress free strain (strain due to growth), see e.g. \cite{Gela,G, 
  Lima0, Lima}.  We assume that the evolution of the nutrient can be described by a reaction-diffusion equation \eqref{sigma}, where $S(\varphi, \sigma)$ is a term accounting for sources and sinks in the nutrient density.  One example is
\begin{align}\label{defn:S}
S(\varphi, \sigma) = - \lambda_c h(\varphi) \sigma + B(\sigma_c - \sigma)
\end{align}
for some non-negative and bounded function $h(\varphi)$, and the coefficient $\lambda_c \geq 0$ has the meaning of a consumption rate.  Meanwhile, the term $B (\sigma_c - \sigma)$ models the supply of nutrients from nearby capillaries, so that $B \geq 0$ is a constant supply rate, and $\sigma_c$ is the nutrient concentration from the capillaries. Furthermore, after non-dimensionalization, the prefactor $\beta > 0$ in front of the time derivative can be interpreted as the ratio between the nutrient diffusion timescale and the tumour doubling timescale.  In many instances, $\beta$ is small, and it makes sense to consider $\beta = 0$ to obtain a quasi-static approximation.  In fact, the same is done for the equation of mechanical stress, which we assume there is an instantaneous relaxation into mechanical equilibrium, leading to the equation \eqref{W,E}.

For boundary conditions, we consider the no-flux condition $\pdnu \mu = 0$, and $\pdnu \varphi = 0$ for the Cahn--Hilliard component.  For $\kappa > 0$, we have a Robin boundary condition for $\sigma$, where $\sigma_B$ can be seen as the nutrient concentration supplied on the boundary, and in the case $\kappa = 0$ we return to a no-flux condition for the nutrient.  For the deformation $\u$, we postulate a zero Dirichlet condition 
on $\GD$ to take into account the possible presence of 
a rigid part of the body such as a bone which prevents variations of the displacement, and a Neumann condition on $\GN$ so that the normal component of the stress $\bsi = W_{,\E}$ on $\GN$ is equal to some load given by a fixed source $\g$.

Let us comment that the dependence of the mobility $m$ on the
displacement $\u$ and the stress $\E(\u)$ is based on the observation that
a tumour induces significant mechanical stresses on the surrounding
tissue during growth, and thus can lead to an inhibition of further
growth \cite{Hormuth}.  A possible volume free energy of the system is
the following:
\begin{align}
E(\varphi, \sigma, \u) &:= \int_\Omega \frac{\eps}{2} \abs{\nx
                             \varphi}^2 + \frac{1}{\eps} \psi(\varphi)
                             + \frac{\beta}{2} \abs{\sigma}^2 +
                             W(\varphi, \E(\u)) d \x, \label{E} \\
 \notag W(\varphi, \E(\u)) & = \frac{1}{2} (\E(\u) - \bE(\varphi)) : \C(\varphi) (\E(\u) - \bE(\varphi)),
\end{align}
where the first two terms of $E$ yields the Ginzburg--Landau energy, the third term is a nutrient free energy, and the last term is the elastic energy.  We note that in the case $\beta = 0$, the nutrient evolves quasi-statically and thus the nutrient free energy is not present in \eqref{E}.

For more information on modelling for tumour growth models in the context of
 Cahn--Hilliard
type models we refer to the book of Cristini and Lowengrub
\cite{CristiniLowengrubbook} and to the articles
\cite{GarckeLNS, GarckeLSS, OdenHawkins}.
Analytical aspects for tumour models based on a
Cahn--Hilliard equation coupled to a nutrient reaction-diffusion equation
have been studied in \cite{ColliGilardiHilhorst, ColliGRSprekels,
  FrigeriGR, GarckeLam1, MiranvilleRS}. Well-posedness result for
extended models of Cahn--Hilliard systems coupled to a flow field have
been studied in \cite{EbenbeckGarcke1, FritzWohlmuthOden, GL,
  LowengrubTiti}. Pioneering numerical simulations showing in particular that
Cahn--Hilliard type models can describe the invasive behaviour of
tumours are due
to Cristini, Lowengrub, Wise and coworkers
\cite{CristiniLowengrub, CristiniLiLowengrubWise, WiseLFC}. For recent numerical
computations for extended models we refer to \cite{agosti, GarckeLNS, oden2}.
Cahn--Hilliard models with mechanical
effects have been first introduced by Cahn and Larch{\'e} \cite{LC} and Onuki
\cite{onuki} and later derived systematically  from
thermodynamic principles by Gurtin \cite{G}. Analytical results for
this so-called Cahn--Larch{\'e} system are due to \cite{BCDGSS, CMA,Gela, GelaSing}
and for a numerical treatment of Cahn--Larch{\'e} systems we refer to
\cite{GW, GraserKS, LLJ}.

Following this Introduction, Section 2 gives an account of the main results
of this paper (existence, uniqueness, regularity and continuous dependence
results). In Section 3 we prove the existence result with a non-standard
Galerkin approach. Section 4 studies a quasi-static limit and in Section 5
further regularity is shown. Under the assumption that the elasticity
tensor is constant we show continuous dependence and uniqueness results in the
final Section 6. 
\section{Main results}
We denote the standard Lebesgue and Sobolev spaces over
$\Omega$ by $L^p := L^p(\Omega)$ and $W^{k,p} := W^{k,p}(\Omega)$, for
$p \in [1,\infty]$ and $k >0$, and denote the corresponding norms by 
$\no{\cdot}_{L^p}$ and $\no{\cdot}_{W^{k,p}}$.  In the case $p = 2$,
we use the notation $H^k := H^k(\Omega) = W^{k,2}(\Omega)$ and the
norm $\no{\cdot}_{H^k}$.  For any Banach space $Z$, we denote its dual
by $Z'$, and the corresponding duality pairing by
$\inn{\cdot}{\cdot}_Z$.  When $Z = H^1(\Omega)$, we use the notation
$\inn{\cdot}{\cdot} = \inn{\cdot}{\cdot}_{H^1}$.  The
$L^2(\Omega)$-inner product is denoted by $(\cdot,\cdot)$, while the
$L^2(\Gamma)$ and $L^2(\GN)$-inner products are denoted by
$(\cdot,\cdot)_{\Gamma}$ and $(\cdot,\cdot)_{\GN}$, respectively.  For
Lebesgue spaces and Sobolev spaces over $\Gamma$, we use the notation
$L^p_\Gamma := L^p(\Gamma)$ and $W^{k,p}_\Gamma := W^{k,p}(\Gamma)$,
respectively, along with the norms $\no{\cdot}_{L^p_\Gamma}$ and
$\no{\cdot}_{W^{k,p}_\Gamma}$.  We define the Sobolev space
$H^2_{\mathbf{n}}(\Omega)$ as the set
$\{f \in H^2(\Omega) \, : \, \pdnu f = 0 \text{ on } \Gamma\}$,
and for the displacement $\u$, we introduce the following function
space:
\begin{align*}
X (\Omega) & := \{ \mathbf{f} \in H^1(\Omega)^d \, : \, \mathbf{f} {\vrule height 5pt depth 4pt\,}_{\GD} = \0 \}.
\end{align*}
Notice that by \cite[Thm.~6.15-4, pp.~409--410]{C}, a Korn-type inequality is valid in $X(\Omega)$: there exists a constant $C_K > 0$ such that 
\begin{align}\label{Korn}
\no{\u}_{H^1} \leq C_K \no{\E(\u)}_{L^2} \quad \forall\, \u \in X(\Omega).
\end{align}
A weak solution to \eqref{Lima} is defined as follows:
\begin{defn}[Weak solution]\label{defn:weak}
We say that $(\varphi, \mu, \sigma, \u)$ is a weak solution to \eqref{Lima} if 
\begin{align*}
\varphi, \sigma & \in L^{2}(0,T;H^1(\Omega)) \cap H^1(0,T;H^1(\Omega)'), \quad \mu \in L^2(0,T;H^1(\Omega)), \\
\u & \in L^2(0,T;\X ),
\end{align*}
with $\varphi(0) = \varphi_0$, $\sigma(0) = \sigma_0$ in $L^2(\Omega)$ and
\begin{subequations}
\begin{alignat}{2}
0 & = \int_0^T \inn{\varphi_t}{\zeta} + (m(\varphi, \sigma, \u, \E( \u)) \nx \mu, \nx \zeta) - (U(\varphi, \sigma, \E(\u)), \zeta) \dt, \label{w:1} \\
 0 &= \int_0^T (\mu, \xi) - \eps (\nx \varphi, \nx \xi) - \eps^{-1} (\psi'(\varphi), \xi) + \chi (\sigma, \xi) - (W_{,\varphi}(\varphi, \E(\u)),  \xi) \dt, \label{w:2} \\
0 & = \int_0^T \beta \inn{\sigma_t}{\zeta} + (\nx \sigma, \nx \zeta) + \kappa(\sigma - \sigma_B,\zeta)_{\Gamma} - (S(\varphi,\sigma), \zeta) \dt, \label{w:3} \\
0 & = \int_0^T (\C(\varphi) (\E(\u) - \bE(\varphi)), \nx \bet) 
%-a( \SS , \nx \bet ) \dt 
-( \g , \bet )_{\GN} \dt 
\label{w:4}
\end{alignat}
\end{subequations}
for all $\zeta \in L^2(0,T;H^1(\Omega))$, $\xi \in L^2(0,T;H^1(\Omega)) \cap L^{\infty}(Q)$ and $\bet \in L^2(0,T; \X)$.
\end{defn}
The first main result of this work concerns the existence of weak solutions to \eqref{Lima} and is formulated as follows.

\begin{thm}[Existence]
\label{thm:exist}
Let $\Omega \subset \R^d$, $d = 2,3$ be a bounded domain with either has 
a $C^{1,1}$-boundary or $\Omega$ is convex.  Suppose 
\begin{enumerate}[label=$(\mathrm{A \arabic*})$, ref = $\mathrm{A \arabic*}$]
\item \label{ass:beta} $\g\in L^2(\GN)^d$, 
%$\s_B \in L^2(\Sigma)$, 
$\eps$ is a positive constant and
$\beta, B, \kappa, \chi$ are non-negative constants such that at least one of $\{B, \kappa\}$ is non-zero if $\beta = 0$.
\item \label{ass:pot} $\psi = \psi_1 + \psi_2$ is
non-negative with $\psi_1, \psi_2 \in C^2(\R)$, $\psi_1$ is convex, $\psi'(0) = 0$, and
\begin{align*}
\forall\, \delta > 0 \, & \exists\, C_\delta > 0 \, : \,
                          \abs{\psi_1'(s)} \leq \delta \psi_1(s) +
                          C_\delta \mbox{ for all } s\in \R, \\
& \exists\, C_1 > 0 \, : \, \abs{\psi_2''(s)} \leq C_1 \mbox{ for all
  } s\in \R.
\end{align*}

\item \label{ass:m} $m \in C^0(\R, \R, \R^{d}, \R^{d\times d})$,
  $f,h, k \in C^0(\R)$ and there exist positive constants
  $C_2, C_3$ such that 
\begin{align*}
0 \leq h(s) \leq 1, \quad \abs{f(s)} \leq 1, \quad \abs{k(s)} \leq 1, 
\quad C_2 \leq m(s, t, \x, \mathcal{A}) \leq C_3
\end{align*}
for all $s, t \in \R$, $\x \in \R^d$, and $\mathcal{A} \in \R^{d \times d}$.

\item \label{ass:W} For the elastic energy $W$, we postulate
\begin{align}
W(s, \E) & = \frac{1}{2} (\E - \bE(s)) : \C(s) (\E - \bE(s)), \quad
           s\in\R, \E \in \R^{d\times d}
\label{def_W}
\end{align}
where the elasticity tensor $\C(s)$ is bounded, Lipschitz continuous and differentiable, while the stress-free 
strain $\bE(s)$ is Lipschitz continuous and differentiable. In
addition, we require that $\C(s)$ fulfills the usual symmetry
conditions of linear elasticity and that there exists a $C_4 >0$
such that for all $s\in\R$, $\E \in \R^{d\times d}_{\rm sym}$
\begin{equation*}
  C_4|\E|^2\le \E : \C (s)\E.
  \end{equation*}
  \item $\lambda_p, \lambda_a, \lambda_c : [0,T] \to \R$ are non-negative continuous functions, while $\sigma_c : \Omega \times (0,T) \to \R$ and $\sigma_B : \Gamma \times (0,T) \to \R$ are non-negative 
  measurable and bounded functions.
\item \label{ass:ini} $\varphi_0 \in H^1(\Omega)$ with $\psi(\varphi_0) \in L^1(\Omega)$ and $\sigma_0 \in L^2(\Omega)$ with 
\begin{align*}
0 \leq \sigma_0 \leq \max \Big ( \no{\sigma_c}_{L^{\infty}(Q)}, \no{\sigma_B}_{L^{\infty}(\Sigma)} \Big) =: M.
\end{align*}
\end{enumerate}
Then, there exists at least one weak solution to \eqref{Lima} in the sense of Definition \ref{defn:weak} such that 
\begin{align}\label{sigma:bdd}
0 \leq \sigma(\x,t) \leq M \quad \forall\, (\x,t) \in Q,
\end{align}
and there exists a positive constant $C = C(\|\varphi_0\|_{H^1},
\|\psi(\varphi_0)\|_{L^1}, \|\g\|_{L^2(\Gamma_N)})$, but not on
$\|\sigma_0\|_{L^2}$ such that 
for a.e.~$s \in (0,T)$, the following energy inequality holds
%\todo{
%\begin{equation}\label{En:ineq}
%\begin{aligned}
%& \no{\nabla  \varphi(s)}_{L^2}^2 + \no{\psi( \varphi(s))}_{L^1} + \beta \no{ \sigma(s)}_{L^2}^2 \\
%& \qquad + \int_\Omega W( \varphi(s), \E( \u(s)) \dx  - \int_{\GN} \g \cdot \u (s) \dA
%+ \int_0^s \no{\nabla  \mu}_{L^2}^2 + \no{ \sigma}_{H^1}^2 \dt
%\\
%& \leq C \int_0^s \abs{\mean{ \mu}} - \no{ \sigma}_{L^2}^2 + \no{ \varphi}_{L^2}^2 \dt \\
%& \qquad + C \Big (1 + \no{\varphi_{0}}_{H^1}^2 + \no{\psi(\varphi_{0})}_{L^1} + \beta \no{\sigma_{0}}_{L^2}^2 \Big ).
%\end{aligned}
%\end{equation}
%[maybe the following fits better]}
\begin{equation}
\begin{aligned}\label{En:ineq}
& \sup_{t \in (0,T)} \Big ( \no{\varphi(t)}_{H^1}^2 
+ \no{\psi(\varphi(t))}_{L^1} + \beta \no{\sigma(t)}_{L^2}^2 + \no{\u(t)}_{\X}^2 \Big ) \\
& \quad 
+ \no{\ph}_{H^1(0,T;H^1(\Omega)')}^2 
+ \beta\no{\s}_{H^1(0,T;H^1(\Omega)')}^2 
+ \no{\mu}_{L^2(0,T;H^1)}^2 
\\ & \quad
+ \no{\sigma}_{L^2(0,T;H^1)}^2 \leq C(1 + \beta\|\sigma_0\|^2_{L^2}).
\end{aligned}
\end{equation}
\end{thm}

It is worth noting that assumption $\eqref{ass:W}$ implies that
$W \in C^1(\R \times \R^{d \times d}, \R)$ is non-negative with
$W(s, \M) = W(s, \M^t)$ for all $s \in \R$ and
$\M \in \R^{d \times d}$, and there exist positive constants $C_4, C_5$
such that for all $s \in \R$, $\M_1, \M_2 \in \R^{d \times d}$ and
$\E \in \R_{\rm sym}^{d \times d}$,
\begin{align}
\label{prop:E:1}
 (W_{,\E}(s, \M_1) - W_{, \E}(s, \M_2)):(\M_1 - \M_2) & \geq C_4 \abs{\M_1 - \M_2}^2, \\
 \label{prop:E:2}
\abs{W(c, \E)} + \abs{W_{,\varphi}(c,\E)} & \leq C_5 ( 1 + \abs{c}^2 + \abs{\E}^2), \\
\label{prop:E:3}
 \abs{W_{,\E}(c,\E)} & \leq C_5  (1 + \abs{c} + \abs{\E}).
\end{align}
Moreover, assumption \eqref{ass:pot} postulates that the derivative 
of the convex part $\psi_1$ can be bounded by $\psi$. This requirement covers the case of regular and polynomial growth potentials, so, for instance, the standard choice $\psi(s)=(s^2-1)^2$ is allowed.

\begin{remark}
Let us remark that the requirement $\u=\0$ on $\SD$ is just to avoid non-necessary technicalities. In fact, the same procedure presented here will be enough to handle the case $\u=\bold f$ for some given source $\bold  f \not =\0$. 
Indeed, it suffices to set $\bold w:= \u - \bold f$ and solve the problem for this auxiliary variable $\bold w$ which now enjoys the same condition as \eqref{SD}.
Moreover, let us claim that also the case 
in which $|\SD|=0$ can be handled by arguing as in \cite{GHab, Gela}.
\end{remark}

It will turn out that the estimates for the solutions are uniform in $\beta \in (0,1]$. 
This allows us to deduce the quasi-static limit $\beta \to 0$ which is formulated as follows.
\begin{thm}[Quasi-static limit]\label{thm:quasi}
For each $\beta \in (0,1]$, let $(\varphi_\beta, \mu_\beta, \sigma_\beta, \u_\beta)$ denote a weak solution to \eqref{Lima} obtained from Theorem \ref{thm:exist} with corresponding initial data $(\varphi_0, \sigma_0)$.  Then, there exist limit functions $(\varphi_*,\mu_*, \sigma_*, \u_*)$ such that, along a non-relabelled subsequence,
\begin{align}
\notag \varphi_\beta & \to \varphi_* \text{ weakly* in } L^{\infty}(0,T;H^1) \cap H^1(0,T;(H^1)') \\
\notag & \hspace{1.2cm} \hbox{and strongly in } C^0([0,T];L^r) \text{ and a.e. in } Q, \\
\notag \mu_\beta & \to  \mu_* \text{ weakly in } L^2(0,T;H^1), \\
\notag \u_\beta & \to \u_* \text{ weakly* in } L^{\infty}(0,T;\X), \\
\notag \notag & \hspace{1.2cm} \hbox{and strongly in } L^2(0,T;\X) \text{ and a.e. in } Q, \\
\notag \beta \sigma_{\beta,t} & \to 0 \text{ weakly in } L^2(0,T;H^1(\Omega)'), \\
\label{sigma:beta} \sigma_\beta & \to \sigma_* \text{ strongly in } L^{2}(0,T;H^1) \text{ and a.e. in } Q,
\end{align}
for any $r < \infty$ in two spatial dimensions and any $r < 6$ in three spatial dimensions.  Furthermore, $(\varphi_*, \mu_*, \sigma_*, \u_*)$ satisfies \eqref{w:1}, \eqref{w:2}, \eqref{w:4} and
\begin{align}\label{qus}
0 & = \int_0^T (\nx \sigma_*, \nx \zeta) + \kappa(\sigma_* - \sigma_B,\zeta)_{\Gamma} - (S(\varphi_*,\sigma_*), \zeta) \dt, 
\end{align}
for all $\zeta \in L^2(0,T;H^1(\Omega))$, $\xi \in L^2(0,T;H^1(\Omega)) \cap L^{\infty}(Q)$ and $\bet \in L^2(0,T;\X)$, the initial condition $\varphi_*(0) = \varphi_0$, the boundedness property \eqref{sigma:bdd}, as well as the energy inequality \eqref{En:ineq} with $\beta = 0$.
\end{thm}

Next, we present a series of regularity assertions for the weak solutions to \eqref{Lima}.

\begin{thm}[Regularity]\label{thm:Reg}
We assume \eqref{ass:beta}-\eqref{ass:ini} and denote by $(\varphi,
\mu, \sigma, \u)$ a weak solution to \eqref{Lima} obtained from Theorem \ref{thm:exist}.  Then, there exists an exponent $p > 2$ such that 
\begin{align}\label{u:Lp}
\nabla \u \in L^\infty(0,T;L^p(\Omega)^{d \times d}).
\end{align}
Moreover, we have the following:
\begin{enumerate}
%\item If $\sigma_0 \in H^1(\Omega)$ and $\sigma_B \in L^2(0,T;H^{1/2}(\Gamma)) \cap H^1(0,T;L^2(\Gamma))$, then 
%\begin{align*}
%\sigma & \in L^\infty(0,T;H^1(\Omega)) \cap L^2(0,T;H^2(\Omega)) \quad \forall\, \beta \in [0,\infty), \\
%\sigma & \in H^1(0,T;L^2(\Omega)) \quad \forall\, \beta > 0.
%\end{align*}
\item If $\sigma_0 \in H^1(\Omega)$ and $\sigma_B \in H^1(0,T;L^2(\Gamma))$, then 
\begin{align*}
\sigma & \in L^\infty(0,T;H^1(\Omega))  \quad \forall\, \beta \in [0,\infty), \quad 
\sigma \in H^1(0,T;L^2(\Omega)) \quad \forall\, \beta > 0.
\end{align*}
Furthermore, if $\Omega$ has a $C^{1,1}$-boundary and $\sigma_B$ also belongs to $L^2(0,T;H^{1/2}(\Gamma))$, then 
\begin{align*}
\sigma & \in L^2(0,T;H^2(\Omega)) \quad \forall\, \beta \in [0,\infty).
\end{align*}
\item Suppose
\begin{enumerate}[label=$(\mathrm{B \arabic*})$, ref = $\mathrm{B \arabic*}$]
\item \label{Vegard} The stress-free strain $\bar \E(\varphi)$ satisfies the affine linear ansatz (Vegard's law)
\begin{align*}
\bar \E(\varphi) = \hat{\E} + \E^* \varphi,
\end{align*}
where $\hat \E$ and $\E^*$ are constant symmetric tensors.
\item \label{C:cons} The elasticity tensor $\C(\ph) = \C$ is a constant, positive definite and symmetric tensor,
\end{enumerate}
hold, then
\begin{align*}
\psi_1'(\varphi)  \in L^2(Q), \quad \varphi \in L^2(0,T;H^2_\bold{n}(\Omega)).
\end{align*}
\end{enumerate}
\end{thm}
\begin{remark}
The regularity assertion \eqref{u:Lp} on the displacement field $\u$ follows from the proof of \cite[Theorem 1.1]{Shi} by choosing $A_{ijkl} = \C(\ph)_{ijkl}$, $f_i = 0$, $f_{ij} = \C(\ph)_{ijkl} (\E(\ph))_{kl}$ and $\tau = g + \C(\ph) \E(\ph)\bm{n}$.  Hence, we omit the details of the proof.
\end{remark}

%\andreaold{
%\begin{remark}
%Let us claim that, following the same lines of \cite{GHab}, 
%it is possible to employ a Caccioppoli argument to slightly improve
%the regularity of $\u$. Indeed we are considering the same elliptic
%equation \eqref{W,E} with a slightly different boundary condition.
%The result would read as follows: assume that $\ph \in L^\gamma(Q)$
%for some $\gamma>2$. Then, there exists some $p\in (2, \gamma]$ 
%such that $\u \in \L\infty {W^{1,p}}$.
%\end{remark}}

Our last results state the continuous dependence of the weak solutions to \eqref{Lima} on the initial conditions and the data, and subsequently leads to the uniqueness of solutions.

\begin{thm}[Continuous dependence]
Further to \eqref{ass:beta}-\eqref{ass:ini},
\eqref{Vegard}-\eqref{C:cons}, we assume 
\begin{enumerate}[label=$(\mathrm{C \arabic*})$, ref = $\mathrm{C \arabic*}$]
\item \label{mob} The mobility $m(\varphi, \sigma, \u, \E(\u))$ is taken to be a constant (w.l.o.g. we set it to be 1).
\item \label{Lip} The functions $f$, $h$ and $k$ are Lipschitz continuous, whose Lipschitz constants we shall denote by a common notation $L > 0$.
\item \label{Cts:pot:generalcase} The convex part $\psi_1$ of the potential $\psi$ satisfies
\begin{align*}
|\psi_1'(s) - \psi_1'(r) | \leq C ( 1+ |s|^q +|r|^q )|s-r| \quad \text{ for all }s,r \in \R,
\end{align*}
and the derivative of the non-convex part $\psi_2$ is Lipschitz continuous (again we denote the Lipschitz constant by $L$).  The exponent $q \in \{2,4\}$ is specified below depending on the norms involved.
\end{enumerate}
Let
\begin{align*}
\mathbb{X}_q := \begin{cases}
L^2  & \text{ if } q = 2, \\
H^1(\Omega)' & \text{ if } q = 4,
\end{cases} \quad \quad \mathcal{P}_q := \begin{cases}
L^\infty & \text{ if } q = 2, \\
L^2 & \text{ if } q = 4.
\end{cases}
\end{align*}
Then, for any pair $\{(\ph_i, \mu_i, \s_i, \u_i)\}_{i=1,2}$ of weak solutions to \eqref{Lima} corresponding to data $\{(\ph_{0,i}, \s_{0,i}, \g_i, \s_{c,i},\s_{B,i})\}_{i=1,2}$, there exists a positive constant $K_1$ not depending on the differences $\ph_1 - \ph_2$, $\mu_1 - \mu_2$, $\s_1 - \s_2$ and $\u_1 - \u_2$, as well as $\beta$, such that
\begin{equation}\label{cts}
\begin{aligned}
&\| \ph_1 - \ph_2\|_{L^\infty(0,T;\mathbb{X}_q) \cap \L2 {H^1}}^2 + \beta \| \s_1 - \s_2\|_{L^\infty(0,T;L^2)}^2 \\
& \qquad + \no{\s_1 - \s_2}_{L^2(0,T;H^1)}^2 +\| \u_1 - \u_2\|_{\mathcal{P}_q(0,T; \X)}^2 + \no{\mu_1 - \mu_2}_{L^2(0,T;\mathbb{X}_q)}^2 \\ 
& \quad \leq K_1 \Big (\| \ph_{0,1} - \ph_{0,2} \|_{\mathbb{X}_q}^2 +\beta \| \s_{0,1} - \s_{0,2} \|_{L^2}^2 \Big ) \\
& \qquad + K_1\Big ( \no{\g_1 - \g_2}_{L^2(\GN)}^2 + \no{\s_{B,1} - \s_{B,2}}_{L^2(\Sigma)}^2 + \no{\s_{c,1} - \s_{c,2}}_{L^2(Q)}^2 \Big ).
\end{aligned}
\end{equation}
In particular, the weak solutions to both \eqref{Lima} and its quasistatic variant are unique.
\end{thm}

Under further assumptions on the convex part $\psi_1$, we obtain the following continuous dependence in stronger norms with a time discretisation approach.

\begin{thm}
Further to \eqref{ass:beta}-\eqref{ass:ini},
\eqref{Vegard}-\eqref{C:cons}, \eqref{mob}-\eqref{Cts:pot:generalcase}
with exponent $q = 2$, we assume that $\Omega$ has a $C^{1,1}$-boundary and 
\begin{enumerate}[label=$(\mathrm{C \arabic*})$, ref = $\mathrm{C \arabic*}$]
\setcounter{enumi}{3}
\item \label{Cts:last} The convex part $\psi_1$ of the potential $\psi$ satisfies 
	\begin{align*}
	\abs{\psi_1''(s) - \psi_1''(r)} \leq C (1 + \abs{s} + \abs{r}) 
	\abs{s-r}, \quad \text{ for all } s, r \in \R.
	\end{align*}
\end{enumerate}
Then, for any pair $\{(\ph_i, \mu_i, \s_i, \u_i)\}_{i=1,2}$ of weak solutions to \eqref{Lima} corresponding to data $\{(\ph_{0,i}, \s_{0,i}, \g_i, \s_{c,i},\s_{B,i})\}_{i=1,2}$, there exists a positive constant $K_2$ not depending on the differences $\ph_1 - \ph_2$, $\mu_1 - \mu_2$, $\s_1 - \s_2$ and $\u_1 - \u_2$, as well as $\beta$, such that
\begin{equation*}
\begin{aligned}
	& \no{\ph_1 - \ph_2}_{L^\infty(0,T;H^1) \cap L^2(0,T;H^2)}^2 + \no{\mu_1 - \mu_2}_{\L2 {H^1}}^2 
	\\
	& \quad \leq 
	K_2 \Big (\| \ph_{0,1} - \ph_{0,2} \|_{H^1}^2 
	+\beta \| \s_{0,1} - \s_{0,2} \|_{L^2}^2 \Big )
	\\ & \qquad 
	+K_2 \Big ( \no{\g_1 - \g_2}_{L^2(\GN)}^2 + \no{\s_{B,1} - \s_{B,2}}_{L^2(\Sigma)}^2
	+ \no{\s_{c,1} - \s_{c,2}}_{L^2(Q)}^2 \Big ).
\end{aligned}
\end{equation*}
Moreover, for any $\beta>0$ and data $\{\s_{0,i}, \s_{B,i}\}_{i=1,2}$ satisfying $\s_{0,i} \in H^1(\Omega)$ and $\s_{B,i} \in H^1(0,T;L^2(\Gamma))$, there is a constant $K_3$ independent of the differences  $\ph_1 - \ph_2$, $\mu_1 - \mu_2$, $\s_1 - \s_2$ and $\u_1 - \u_2$, such that
\begin{align*}
	 \no{\sigma_1 - \s_2}_{\H1 {L^2} \cap L^\infty(0,T;H^1) %\cap L^2(0,T;H^2)
	}^2
%	\\
%	& \quad 
& \leq 
	K_3 \Big (\| \ph_{0,1} - \ph_{0,2} \|_{H^1}^2 
	+ \| \s_{0,1} - \s_{0,2} \|_{H^1}^2 + \no{\g_1 - \g_2}_{L^2(\GN)}^2 \Big )
	\\ & \qquad 
	+K_3 \Big ( \no{\s_{B,1} - \s_{B,2}}_{H^1(0,T;L^2_\Gamma)}
	 %\cap L^2(0,T;H^{1/2}(\Gamma))}^2
	+ \no{\s_{c,1} - \s_{c,2}}_{L^2(Q)}^2 \Big ).
\end{align*}
Lastly, if $\s_{B,1},\s_{B,2}$ also belong to $L^2(0,T;H^{1/2}(\Gamma))$, we also have
\begin{align*}
	 \no{\sigma_1 - \s_2}_{L^2(0,T;H^2)}^2 
	%\\ & \quad 	
	& \leq 
	\andrea{K_3} \Big (\| \ph_{0,1} - \ph_{0,2} \|_{H^1}^2 
	+ \| \s_{0,1} - \s_{0,2} \|_{H^1}^2 + \no{\g_1 - \g_2}_{L^2(\GN)}^2 \Big )
	\\ & \qquad 
	+\andrea{K_3} \Big ( \no{\s_{B,1} - \s_{B,2}}_{H^1(0,T;\andrea{L^2_\Gamma)}
	 \cap L^2(0,T;H^{1/2}(\Gamma))
	 }^2
	+ \no{\s_{c,1} - \s_{c,2}}_{L^2(Q)}^2 \Big ).
\end{align*}
\end{thm}
Notice that conditions \eqref{Cts:pot:generalcase} and \eqref{Cts:last} still comply for the 
classical quartic potential $\psi(s) = (s^2 - 1)^2$.

\section{Existence}
Due to the presence of the source terms $U$ and $S$, the system \eqref{Lima} does not admit a variational structure, and so an implicit time discretisation such as the one used in \cite{GHab,Gela} may no longer be applicable.  Hence, we consider a Faedo--Galerkin approximation to establish a weak solution to system \eqref{Lima}. 
\subsection{Galerkin approximation}
Let us point out that, since \eqref{ass:pot} allows for the potential $\psi$ to have arbitrary polynomial growth, in the first step we prove Theorem \ref{thm:exist} with a more regular initial condition $\varphi_0 \in H^2_{\bold n}(\Omega)$, and then in Section \ref{sec:ini} show how to complete the proof for an initial condition satisfying just \eqref{ass:ini}.  To this end, let us consider
\begin{itemize}
\item $\{z_i\}_{i \in \N}$ as the set of eigenfunctions of the Neumann-Laplacian operator that is orthonormal in $L^2(\Omega)$ and orthogonal in $H^1(\Omega)$ with $z_1$ is the constant function $(\frac 1{|\Omega|})^{1/2}$ and $(z_i, 1) = 0$ for $i \geq 2$. In \cite[\S 3]{GL} it is also shown that $\{z_i\}_{i \in \N}$ forms a basis of $H^2_{\bold n}(\Omega)$;
\item $\{\y_i\}_{i \in \N}$ as a Schauder basis of $\X$, {see \cite{Alt}. One can choose for example eigenfunctions of a corresponding boundary value problem for an elasticity system which leads to a basis orthogonal in $L^2(\Omega)^d$ (see \cite[Thm.~3.12.1, pp.~219-220]{LVC}).}
\end{itemize}
Next, we define finite-dimensional spaces $Z_k$ and $\bY_k$
as the linear span of the first $k$ functions of $\{z_i\}_{i \in \N}$
and $\{\y_i\}_{i \in \N}$, respectively, and we denote by $\Pi_k$ the $L^2$-projection onto the space $Z_k$.  Then, the Faedo--Galerkin
approximation of \eqref{w:1}-\eqref{w:4} reads as: for any $k \in \N$
find $(\varphi_k, \mu_k, \sigma_k, \u_k)$ of the form
\begin{align*}
\varphi_k = \sum_{i=1}^k a_{i}^k(t) z_i(\x), \quad \mu_k = \sum_{i=1}^k b_{i}^k(t) z_i(\x), \quad \sigma_k = \sum_{i=1}^k c_{i}^k(t) z_i(\x), \quad \u_k = \sum_{i=1}^k d_{i}^k(t) \y_i(\x), 
\end{align*}
satisfying for a.e. $t \in (0,T)$ and $j \in \{1, \dots, k\}$,
\begin{subequations}\label{Gal:a}
\begin{alignat}{2}
0 & = (\varphi_k', z_j) + (m(\varphi_k, \sigma_k, \u_k, \E(\u_k)) \nx \mu_k, \nx z_j) - (\tilde{U}(\varphi_k,\sigma_k, \E(\u_k)), z_j), \label{G:a:phi} \\
\tilde{U}& = \lambda_p f(\varphi_k) g(\sigma_k) / (1 + \abs{W_{,\E}(\varphi_k, \E(\u_k))}) - \lambda_a k(\varphi_k), \\
\label{G:a:mu} 0 & = (\mu_k - \eps^{-1} \psi'(\varphi_k) + \chi \sigma_k - W_{,\varphi}(\varphi_k, \E(\u_k)), z_j) - (\eps \nx \varphi_k, \nx z_j) \\
0 & =\beta (\sigma_k', z_j) + (\nx \sigma_k, \nx z_j) + \kappa(\sigma_k - \sigma_B, z_j)_\Gamma - (S(\varphi_k, \sigma_k), z_j), \label{G:a:sig} \\
S& = -\lambda_c h(\ph_k)\s_k + B (\s_c-\s_k),  \\
0 & = (W_{,\E}(\varphi_k, \E(\u_k)), \nx \y_j)
- ( \g, \y_j)_{\GN} 
, \label{G:a:u}\\
\varphi_k(0) & = \varphi_{k,0} := \Pi_k \varphi_{0}, 
\quad 
\sigma_{k}(0) = \sigma_{k,0} := \Pi_k \sigma_0,
\end{alignat}
\end{subequations}
where in the definition of $\tilde U$, the function
$g(s) = \max(0, \min(s,\no{\sigma_B}_{L^{\infty}(\Sigma}),
\no{\sigma_c}_{L^{\infty}(Q)}))$ is a truncation.  It will turn out
that the nutrient equation satisfies a comparison principle, but this
is not valid at the Galerkin level, and thus we introduce the
truncation $g$ to first derive the necessary a priori estimates, and
then remove it at the continuous level.

The orthogonality of $\{z_i\}_{i \in \N}$ 
with respect to the $L^2$-inner product allows us to
express \eqref{Gal:a} as a system of ordinary differential equations
in the coefficient vectors $\bm{a} := (a_1^k, \dots, a_k^k)$,
$\bm{b} := (b_1^k, \dots, b_k^k)$,
$\bm{c} := (c_1^k, \dots, c_k^k)$ and
$\bm{d}: = (d_1^k, \dots, d_k^k)$.  It is not hard to see
that the continuity of $m$, $\psi'$, $f$, $g$, $h$, $k$,
$W_{,\varphi}$, $W_{,\E}$ with respect to their arguments, as well
as the continuity of $\lambda_p$, $\lambda_a$ and $\lambda_c$ with
respect to time yield that the differential-algebraic system contains only contributions that are continuous in
$\bm{a}, \bm{b}, \bm{c}, \bm{d}$.

Here, we are going to show that the above system can be
expressed as a system of ODEs in terms of $\bm a$ and $\bm c$ only
and that classical results ensure the existence of a regular
solution.  From an analysis of equation \eqref{G:a:mu} it is
straightforward to realise that $\bm b$ can be expressed as a function
of $\bm{a}, \bm{c}$, and $ \bm{d}$ in a continuously differentiable
fashion.  Moreover, let us claim that $\bm{d}$ can be
expressed as a function of $\bm a$ only.  In the direction of
formalising this fact, let us fix for convenience the following
notation:
\begin{align*}
	\ph(\bm a) := \ph_k = \sum_{i=1}^k a_{i}^k z_i(\x), 
	\quad 
	\u (\bm d) := \u_k = \sum_{i=1}^k d_{i}^k \y_i(\x).
\end{align*}
Furthermore, we point out that equation \eqref{G:a:u} can be written as $\0 = F(\bm a, \bm d),$ for a function $F:\R^k \times\R^k \to \R^k $ defined by 
\begin{align*}
	F(\bm a, \bm d)_{i} &= (W_{,\E}(\ph(\bm a), \E(\u(\bm d))), \nx \y_i)
	- (\g, \y_i)_{\GN}
	\\
	&= \big(\sum_{j=1}^{k} \C(\ph(\bm a))\E(d_j^k \y_j), \nx \y_i \big)
	- \big(\sum_{j=1}^{k} \C(\ph(\bm a))\bar\E(a_j^k z_j), \nx \y_i \big)
	\\
	& \quad
	- (\g, \y_i)_{\GN}, \quad \hbox{for $i=1,...,k$}.
\end{align*}
Moreover, by virtue of symmetry, we can replace $\nx \y_i$ with $\E(\y_i)$ to infer that 
\begin{align*}
	F(\bm a, \bm d)_{i} &= 
	\big(\sum_{j=1}^{k} \C(\ph(\bm a))\E(d_j^k \y_j), \E(\y_i) \big)
	- \big(\sum_{j=1}^{k} \C(\ph(\bm a))\bar\E(a_j^k z_j), \E(\y_i) \big)
	- (\g, \y_i)_{\GN}
	\\
	&=: (\mathbb{A}(\bm a) \bm d)_i - \bm q (\bm a)_i,
\end{align*}
where
\begin{align*}
	(\mathbb{A}(\bm a) \bm d)_i &= 
	\big(\sum_{j=1}^{k} \C(\ph(\bm a))\E(d_j^k \y_j), \E(\y_i) \big),
	\\ 
	\bm q (\bm a)_i &=  \big(\sum_{j=1}^{k} \C(\ph(\bm a))\bar\E(a_j^k z_j), \E(\y_i) \big)
	+ (\g, \y_i)_{\GN}.
\end{align*}
An easy calculation using the fact that the tensor $\C$ is positive definite shows that the matrix $\mathbb{A}(\bm a)$ is positive definite and hence invertible so that we can uniquely solve the linear system
\begin{align}
	F(\bm a, \bm d) = \mathbb{A}(\bm a) \bm d - \bm q (\bm a)=\0
	\label{d:1}
\end{align}
and deduce that
\begin{align*}
	\bm d = \mathbb{A}^{-1}(\bm a) \bm q (\bm a),
\end{align*}
which in turn proves that the solution $\bm d$ of \eqref{d:1} depends continuously on $\bm a$.  Recalling that $\bm b$ can be expressed as a function of $\bm a, \bm c,$ and $\bm d$, we can rephrase equations 
\eqref{G:a:phi} and \eqref{G:a:sig} as a system of ordinary differential equations:
\begin{align*}
		\begin{cases}
		\bm {a'} = H_1 (\bm a, \bm c),
		\\
		\bm {c'} = H_2 (\bm a, \bm c),		
		\end{cases}
\end{align*}
for suitable functions $H_i$, $i=1,2$, that are continuous with respect to their arguments. In light of the above observations, we invoke the Cauchy--Peano theorem to obtain the existence of $T_k \in (0,T]$ and local solutions $\bm{a}, \bm{c} \in C^1([0,T_k], \R^k)$ solving \eqref{G:a:phi} and \eqref{G:a:sig}, from which we also get $\bm{b}, \bm{d} \in C^1([0,T_k], \R^k)$.  Moreover, we can endow an initial condition for the Galerkin approximation $\u_k$ via the relation \eqref{d:1}. Namely, we set
\begin{align*}
\bm{d}(0) = \mathbb{A}^{-1}(\bm{a}(0))\bm{q}(\bm{a}(0)) 
\end{align*}
which implies that $\u_{k,0} := \u_k(0) = \sum_{i=1}^kd_i^k(0) \y_i(\x) \in \X$ satisfies
\begin{align}\label{uk:0}
( W_{,\E}(\varphi_{k,0}, \E(\u_{k,0})), \nabla \y_j) = (\g, \y_j)_{\GN} \quad \forall\, 1 \leq j \leq k.
\end{align}
Then, multiplying the above by $d_j^k(0)$, summing from $j = 1$ to $k$, and employing \eqref{prop:E:1}, Korn's inequality and the trace theorem yields
%\old{\begin{align*}
%c \no{\u_{k,0}}_{\X}^2  \leq C_2 \no{\E(\u_{k,0})}_{L^2}^2 & \leq (\C(\varphi_{k,0}) (\E(\u_{k,0}) - \bar \E(\varphi_{k,0})), \E(\u_{k,0})) \\
%& = (\C(\varphi_{k,0}) \bar \E(\varphi_{k,0}), \E(\u_{k,0})) + (\g, \u_{k,0})_{\GN} \\
%& \leq C \Big ( \no{\varphi_{k,0}}_{L^2} + \no{\g}_{L^2(\GN)} \Big ) \no{\u_{k,0}}_{\X},
%\end{align*}}
\begin{align*}
\frac {C_4}{C_K^2} \no{\u_{k,0}}_{\X}^2  \leq C_4 \no{\E(\u_{k,0})}_{L^2}^2 & \leq (\C(\varphi_{k,0}) \E(\u_{k,0}) , \E(\u_{k,0})) \\
& = (\C(\varphi_{k,0}) \bar \E(\varphi_{k,0}), \E(\u_{k,0})) + (\g, \u_{k,0})_{\GN} \\
& \leq C \Big ( \no{\varphi_{k,0}}_{L^2} + \no{\g}_{L^2(\GN)} \Big ) \no{\u_{k,0}}_{\X},
\end{align*}
which implies the following estimate of $\u_{k,0}$ in $\X$:
\begin{align}\label{uk:0:est}
\no{\u_{k,0}}_{\X} \leq C\Big ( 1 + \no{\varphi_{k,0}}_{L^2}  \Big )
\end{align}
for a positive constant $C$ independent of $k$.

Next, our goal is to derive uniform estimates in $k$ to pass to the
limit.   In the sequel, we denote positive constants that are
independent of $\beta$ and $k$, which may vary from line to line, by the symbol $C$. 

\subsection{A priori estimates}\label{sec:apr}
Let $R, K > 0$ be constants yet to be determined.
Multiplying \eqref{G:a:phi} with $b_j^k$ and with $K a_j^k$,
\eqref{G:a:mu} with $(a_j^k)'$, \eqref{G:a:sig} with
$R c_j^k $ and \eqref{G:a:u} with $(d_j^k)'$.  Summing from $j=1$
to $k$ and using the symmetry of the elastic tensor $\C $ yields
\begin{align*}
  0 & = (\varphi_k', \mu_k) + \no{m^{1/2} \nx \mu_k}_{L^2}^2 - (\tilde U, \mu_k), \\
  0 & = \frac{d}{dt} \frac{K}{2} \no{\varphi_k}_{L^2}^2 + K (m \nx \mu_k, \nx \varphi_k) - K (\tilde U, \varphi_k), \\
  0 & = (\mu_k + \chi \sigma_k - W_{,\varphi}(\varphi_k, \E(\u_k)), \varphi_k') - \frac{d}{dt} \int_\Omega \frac{1}{\eps} \psi(\varphi_k) + \frac{\eps}{2} \abs{\nx \varphi_k}^2 
      d\x, \\
  0 & = \frac{d}{dt} \frac{R\beta}{2} \no{\sigma_k}_{L^2}^2 +
R\no{\nx \sigma_k}_{L^2}^2 + R \kappa
      \no{\sigma_k}_{L^2_\Gamma}^2 - R \kappa (\sigma_B,
      \sigma_k)_\Gamma- R(S(\varphi_k, \sigma_k), \sigma_k)\\
  0 & = 
      (W_{,\E}(\varphi_k, \E(\u_k)), \E(\u_k')) 
      - ( \g, \u_k')_{\GN}. 
\end{align*}
Taking note of the identity
\begin{align*}
(W_{,\varphi}, \varphi_k') + (W_{,\E}, \E(\u_k')) = \frac{d}{dt} \int_\Omega W(\varphi_k, \E(\u_k)) \dx,
\end{align*}
we obtain the following after summing the equations
\begin{equation}\label{Energy}
\begin{aligned}
  & \frac{d}{dt} \int_\Omega \frac{\eps}{2} \abs{\nx \varphi_k}^2 +
  \frac{1}{\eps} \psi(\varphi_k) + \frac{K}{2} \abs{\varphi_k}^2 +
  \frac{R\beta}{2} \abs{\sigma_k}^2
  d\x \\
  & \qquad + \frac{d}{dt} \int_\Omega W(\varphi_k, \E(\u_k)) d\x -
  \frac{d}{dt} \int_{\GN} \g \cdot \u_k \dA
  + \no{m^{1/2} \nx \mu_k}_{L^{2}}^2 \\
  & \qquad + R \no{\nx \sigma_k}_{L^2}^2 + R \kappa \no{\sigma_k}_{L^2_\Gamma}^2 + R \lambda_c \no{h^{1/2} \sigma_k}_{L^2}^2 + R B \no{\sigma_k}_{L^{2}}^2 \\
  & \quad = (\tilde U, \mu_k + K \varphi_k) + (\chi \sigma_k,
  \varphi_k') - K(m \nx \mu_k, \nx \varphi_k) + \andrea{R} B(\sigma_c,
  \sigma_k) + \andrea{R} \kappa (\sigma_B, \sigma_k)_\Gamma.
\end{aligned}
\end{equation}
The last two terms on the right-hand side can be handled using the
Cauchy--Schwarz and Young inequalities as follows:
\begin{align*}
B(\sigma_c, \sigma_k) + \kappa (\sigma_B, \sigma_k)_\Gamma \leq 
\tfrac{1}{2} \kappa ( \no{\sigma_k}_{L^2_\Gamma}^2 + \no{\sigma_B}_{L^2_\Gamma}^2) +  d_1 \no{\sigma_k}_{L^2}^2 + \tfrac{1}{4} \tfrac{B^2}{d_1} \no{\sigma_c}_{L^2}^2,
\end{align*}
where $d_1$ is a positive constant to be determined later.  Let us observe that for an arbitrary test function $\zeta \in H^1(\Omega)$, we multiply \eqref{G:a:phi} with the coefficients of $\Pi_k \zeta$, leading to 
\begin{align*}
  (\varphi_k', \zeta) = (\varphi_k', \Pi_k \zeta) = -(m \nx \mu_k, \nx \Pi_k \zeta) + (\tilde U, \Pi_k \zeta).
\end{align*}
The boundedness of $m$ and $\tilde U$, and the estimate $\no{\Pi_k \zeta}_{H^1} \leq C \no{\zeta}_{H^1}$ implies that
\begin{align}\label{varphik'}
\no{\varphi_k'}_{(H^1)'} \leq \andrea{C_3} \no{\nx \mu_k}_{L^{2}} + C_{\lambda,g},
\end{align}
for some positive constant $C_{\lambda,g}$ depending only on
$\max_{t\in [0,T]} \lambda_p(t)$, $\max_{t\in [0,T]} \lambda_a(t)$ and
$\max_{s\in \R} g(s)$.  We now estimate \andrea{the} remaining terms
of the right-hand side as follows:
\begin{align*}
  (\tilde U, \mu_k) & = (\tilde U , \mu_k - \mean{\mu_k}) + \mean{\mu_k} (\tilde U , 1)\\  & 
  \leq C_{\lambda, g} \no{\mu_k - \mean{\mu_k}}_{L^1} + C_{\lambda,g} \abs{\mean{\mu_k}} \\
                    & \leq \eta \no{\nx \mu_k}_{L^2}^2 + C(\eta^{-1}, c_p, C_{\lambda, g}) + C_{\lambda,g} \Big (\abs{\mean{\mu_k}} - \no{\sigma_k}_{L^2}^2 + \no{\sigma_k}_{L^2}^2 \Big), \\
  (\tilde U, K \varphi_k) & \leq C(K, C_{\lambda, g}) + \no{\varphi_k}_{L^2}^2, \\
  K(m \nx \mu_k, \nx \varphi_k) & \leq \eta \no{\nx \mu_k}_{L^2}^2 + C(\eta^{-1}, K, \andrea{C_3}) \no{\nx \varphi_k}_{L^2}^2, \\
  (\sigma_k, \varphi_k') & \leq \no{\sigma_k}_{H^1} \no{\varphi_k'}_{(H^1)'} \leq \eta \no{\nx \mu_k}_{L^2}^2 + C(\eta^{-1}, \andrea{C_3}) + C(\eta^{-1}, C_{\lambda,g})\no{\sigma_k}_{H^1}^2,
\end{align*}
where $c_p > 0$ is the constant from the Poincar\'{e} inequality and
$\eta$ is a constant yet to be determined and where in the last
inequality we made use of \eqref{varphik'}.  Putting everything together and using the lower bound for the mobility, we obtain from \eqref{Energy}
\begin{equation}\label{est1a}
\begin{aligned}
& \frac{d}{dt} \int_\Omega \frac{\eps}{2} \abs{\nx \varphi_k}^2 + \frac{1}{\eps} \psi(\varphi_k) + \frac{K}{2} \abs{\varphi_k}^2 + \frac{R\beta}{2} \abs{\sigma_k}^2 
d\x \\
& \qquad 
+ \frac{d}{dt} \int_\Omega W(\varphi_k, \E(\u_k)) d\x 
- \frac{d}{dt} \int_{\GN} \g \,\cdot \u_k \dA \\
& \qquad 
+ ( \andrea{C_2} -3 \eta) \no{\nx \mu_k}_{L^{2}}^2 + \andrea{R}
\no{\nx \sigma_k}_{L^2}^2 + \frac{\andrea{R\kappa}}{2} \no{\sigma_k}_{L^2_\Gamma}^2  \\
& \quad \leq C(K, \eta^{-1}, C_{\lambda, g}, \kappa, d_1^{-1}, \sigma_B, \sigma_c) + (\andrea{R}d_1 + C(\eta^{-1}, C_{\lambda,g})) \no{\sigma_k}_{H^1}^2 \\
& \qquad  + C_{\lambda,g} \Big (\abs{\mean{\mu_k}} - \no{\sigma_k}_{L^2}^2 \Big) + \no{\varphi_k}_{L^2}^2 + C(\eta^{-1}, K) \no{\nx \varphi_k}_{L^2}^2.
\end{aligned}
\end{equation}
Let us point out that if $K = 0$ then the last term on the right-hand side involving $\no{\nabla \varphi_k}_{L^2}^2$ vanishes.  Furthermore, we recall the generalised Poincar\'{e} inequality: There exists a positive constant $C_p = C_p(\Omega)$ such that for all $f \in H^1(\Omega)$,
\begin{align}\label{Poin}
\no{f}_{H^1}^2 \leq C_p \Big ( \no{\nabla f}_{L^2}^2 + \no{f}_{L^2_\Gamma}^2 \Big ),
\end{align}
so that by choosing $d_1 < \min(1,\kappa) \frac{1}{8C_p}$ and $\andrea{R} > \frac{8 C_p C(\eta^{-1}, C_{\lambda, g})}{\min(1,\kappa)}$, we obtain
\begin{align*}
  (\andrea{R} d_1 + C(\eta^{-1}, C_{\lambda, g}))  \no{\sigma_k}_{H^1}^2 & \leq \frac{\andrea{R}}{4 C_p} \min(1, \kappa) \no{\sigma_k}_{H^1}^2 \leq \frac{\andrea{R}}{4} \no{\nabla \sigma_k}_{L^2}^2 + \frac{\andrea{R\kappa}}{4} \no{\sigma_k}_{L^2_\Gamma}^2.
\end{align*} 
Then, in \eqref{est1a} we choose $\eta = \frac{\andrea{C_2}}{4}$ 
\begin{equation}\label{Energy:ptl}
\begin{aligned}
& \frac{d}{dt} \int_\Omega \frac{\eps}{2} \abs{\nabla \varphi_k}^2 + \frac{1}{\eps} \psi(\varphi_k) + \frac{K}{2} \abs{\varphi_k}^2 + \frac{\andrea{R}\beta}{2} \abs{\sigma_k}^2 
d\x \\
& \qquad 
+ \frac{d}{dt} \int_\Omega W(\varphi_k, \E(\u_k)) d\x 
- \frac{d}{dt} \int_{\GN} \g \andrea{\cdot}\u_k \dA \\
& \qquad 
+ \frac{\andrea{C_2}}{4} \no{\nx \mu_k}_{L^{2}}^2 + \frac{3 \andrea{R}}{4} \no{\nx \sigma_k}_{L^2}^2 + \frac{\andrea{R\kappa}}{4} \no{\sigma_k}_{L^2_\Gamma}^2  \\
& \quad \leq C + C_{\lambda, g} \Big ( \abs{\mean{\mu_k}} - \no{\sigma_k}_{L^2}^2 \Big ) + \no{\varphi_k}_{L^2}^2 + C(\eta^{-1}, K) \no{\nabla \varphi_k}_{L^2}^2.
\end{aligned}
\end{equation}
We will use this expression to obtain an energy inequality for the limit solutions.  Next, to estimate the term involving the mean value $\mean{\mu_k}$, we choose $j = 1$ in \eqref{G:a:mu} (recalling $z_1$ is constant), and use \eqref{ass:beta}-\eqref{ass:m} all together and \eqref{ass:W} to deduce that
\begin{equation}\label{mean:mu}
\begin{aligned}
\abs{\mean{\mu_k}}& \leq C \Big ( 1 + \no{\psi_1'(\varphi_k)}_{L^1} + \no{\varphi_k}_{L^2}^2 + \no{\E(\u_k)}_{L^2}^2 \Big ) + \no{\sigma_k}_{L^2}^2 \\
& \leq C \Big ( 1 + \no{\psi(\varphi_k)}_{L^1} + \no{\varphi_k}_{L^2}^2  + \no{\E(\u_k)}_{L^2}^2 \Big ) + \no{\sigma_k}_{L^2}^2,
\end{aligned}
\end{equation}
and \andrea{so, using \eqref{Poin},} \eqref{Energy:ptl} simplifies to
\begin{equation}\label{est1c}
\begin{aligned}
& \frac{d}{dt} \int_\Omega \abs{\nx \varphi_k}^2 + \psi(\varphi_k) + K \abs{\varphi_k}^2 + \beta \abs{\sigma_k}^2 
d\x \\
& \qquad
+ \frac{d}{dt} \int_\Omega W(\varphi_k, \E(\u_k)) \dx 
- \frac{d}{dt} \int_{\GN} \g \cdot \u_k \dA
+ \no{\nx \mu_k}_{L^{2}}^2 + \no{\sigma_k}_{H^1}^2   \\
& \quad \leq C \Big (1 + \no{\psi(\varphi_k)}_{L^1} + \no{\varphi_k}_{H^1}^2 + \no{\E(\u_k)}_{L^2}^2 \Big ).
\end{aligned}
\end{equation}
Using the strict monotonicity of $W_{,\E}$ with respect to its second argument, we find that for any $s \in \R$ and $\M \in \R^{d \times d}$
\begin{align*}
W(s, \M)  = W(s, \0) + \int_0^1 W_{,\E}(s, t \M) : t\M \frac{1}{t} \dt \geq \frac{C_\andrea{4}}{2} \abs{\M}^2 - C(1 + \abs{s}^2).
\end{align*}
Hence, by Young's inequality, the trace theorem, and Korn's inequality
we have for some positive constants $c_1$ and $c_2$ that
\begin{align}\label{W:lb}
\int_\Omega W(\varphi_k, \E(\u_k)) d\x 
- \int_{\GN} \g \cdot \u_k \dA 
\geq c_1 \no{\E(\u_k)}_{L^2}^2 - c_2( 1 + \no{\varphi_k}_{L^2}^2 ).
\end{align}
On the other hand, by combining \eqref{prop:E:2} with \eqref{uk:0:est}, we deduce
\begin{equation}\label{W:0:est}
\begin{aligned}
& \int_\Omega W(\varphi_{k,0}, \E(\u_{k,0})) d\x - \int_{\GN} \g \cdot \u_{k,0} \dA \\
& \quad \leq C \Big (1 + \no{\varphi_{k,0}}_{L^2}^2 + \no{\u_{k,0}}_{\X}^2 + \no{\g}_{L^2(\GN)}^2 \Big )  \leq C \Big ( 1 + \no{\varphi_{k,0}}_{L^2}^2 \Big ).
\end{aligned}
\end{equation}

Therefore, \andrea{provided we choose} $K > c_2$, integrating \eqref{est1c} in time and employing the above estimate leads to
\begin{equation}\label{est1d}
\begin{aligned}
& \Big (\no{\varphi_k}_{H^1}^2 
+ \no{\psi(\varphi_k)}_{L^1} + \beta \no{\sigma_k}_{L^2}^2 + \no{\E(\u_k)}_{L^2}^2  \Big)(t)  \\
& \qquad
+ \no{\nx \mu_k}_{L^2(0,t;L^2)}^2 + \no{\sigma_k}_{L^2(0,t;H^1)}^2 \\
& \quad \leq C\Big ( 1 + \no{\psi(\varphi_k)}_{L^1(0,t;L^1)} + \no{\varphi_k}_{L^2(0,t;H^1)}^2 + \no{\E(\u_k)}_{L^2(0,t;L^2)}^2 \Big ) \\
& \qquad + C \no{\varphi_{k,0}}_{H^1}^2 
+ \no{\psi(\varphi_{k,0})}_{L^1} + \beta \no{\sigma_{k,0}}_{L^2}^2  \quad \forall\, t \in (0,T).
\end{aligned}
\end{equation}
Thanks to the \andreaold{fact that the} $\{z_i\}_{i \in \N}$
\andreaold{are a basis} in $H^2(\Omega)$ and \andreaold{are
  orthonormal} in $L^2(\Omega)$, and recalling our assumptions on the set $\Omega$,
   there exists a positive constant $C$
such that $\no{\varphi_{k,0}}_{H^2} \leq C \no{\varphi_0}_{H^2}$ and
$\no{\sigma_{k,0}}_{L^2} \leq \no{\sigma_0}_{L^2}$.  Furthermore,
since $\varphi_0 \in H^2_{\mathbf{n}}(\Omega)$, by \cite[\S 3, p.~329]{GL} we
have $\varphi_{k,0} = \Pi_k \varphi_0 \to \varphi_0$ strongly in
$H^2_{\mathbf{n}}(\Omega) \subset L^{\infty}(\Omega)$ and a.e.~in $\Omega$. This
implies that $\varphi_{k,0}$ is bounded uniformly in
$L^{\infty}(\Omega)$, and in turn we can deduce the existence of a
positive constant $c_*$ such that
\begin{align*}
\no{\psi(\varphi_{k,0})}_{L^{\infty}} \leq c_* \text{ for all } k \in \N.
\end{align*}
%Then, given any arbitrary $\eta > 0$, we can find measurable subsets $E \subset \Omega$ with $\abs{E} < \frac{\eta}{c_*}$ such that 
%\begin{align*}
%\int_{E} \psi(\varphi_{k,0}) d\x \leq \no{\psi(\varphi_{k,0})}_{L^{\infty}} \abs{E} < \eta.
%\end{align*}
%This shows that the family $\{\psi(\varphi_{k,0})\}_{k \in \N}$ is uniformly integrable over $\Omega$, and by Vitali's 
Then, by Lebesgue convergence theorem we obtain
\begin{align}\label{psi:ini:conv}
\int_{\Omega} \psi(\varphi_{k,0}) d\x \to \int_{\Omega} \psi(\varphi_0) d\x \quad \text{ as } k \to \infty.
\end{align}
Hence, there exists a positive constant $C$ such that
\begin{align*}
\no{\psi(\varphi_{k,0})}_{L^1} \leq C,
\end{align*}
and through the use of Gronwall's inequality in integral form and Korn's inequality \eqref{Korn} we obtain
\begin{equation}\label{Est1}
\begin{aligned}
& \sup_{t \in (0,T)} \Big ( \no{\varphi_k(t)}_{H^1}^2 
+ \no{\psi(\varphi_k(t))}_{L^1} + \beta \no{\sigma_k(t)}_{L^2}^2 + \no{\u_k(t)}_{\X}^2 \Big ) \\
& \quad + \no{\nx \mu_k}_{L^2(0,T;L^2)}^2 
+ \no{\sigma_k}_{L^2(0,T;H^1)}^2 \leq C(1 + \beta \norma{\sigma_0}^2).
\end{aligned}
\end{equation}
Then, from \eqref{mean:mu} we have $\mean{\mu_k}(t)$ is bounded in $L^{\infty}(0,T)$, which leads to
\begin{align}\label{Est2}
\no{\mu_k}_{L^2(0,T;H^1)} \leq C(1 +\beta \norma{\sigma_0}^2).
\end{align}
Furthermore, from \eqref{varphik'} and arguing similarly with \eqref{G:a:sig}, we also have
\begin{align}\label{Est3}
\no{\varphi_k'}_{L^2(0,T;(H^1)')} + \beta \no{\sigma_k'}_{L^2(0,T;(H^1)')} \leq C(1 + \beta \norma{\sigma_0}^2).
\end{align}

\subsection{Compactness assertions and passing to the limit}\label{sec:compact}
From the above estimates we immediately deduce the existence of functions $(\ph, \mu, \u, \sigma)$ such that, for a non-relabelled subsequence, we have
\begin{equation*}
\begin{aligned}
\ph_k & \to \varphi \text{ weakly* in } L^{\infty}(0,T;H^1) \cap H^1(0,T;(H^1)') \\
%\ph_k & \to \varphi 
 & \hspace{1cm}\andrea{\text{and strongly in } C^0([0,T];L^r) \text{ and a.e. in } Q,} \\
\mu_k & \to \mu \text{ weakly in } L^2(0,T;H^1), \\
\u_k & \to \u \text{ weakly* in } L^{\infty}(0,T;{\X}), \\
\sigma_k & \to \sigma \text{ weakly* in } L^2(0,T;H^1) \cap  L^{\infty}(0,T;L^2) \cap H^1(0,T;(H^1)') \\
%\sigma_k & \to \sigma 
 & \hspace{1cm}\andrea{\text{and strongly in } L^2([0,T];L^r) \text{ and a.e. in } Q}
\end{aligned}
\end{equation*}
for any $r < \infty$ in two spatial dimensions and any $r < 6$ in
three spatial dimensions. Next, we deduce the strong convergence of
$\u_k$ to $\u$, which can be argued as follows.  Since the span of
$\{\y_i\}_{i \in \N}$ is dense in $\X$, we can choose a sequence $\{
\bm{v}_k\}_{k \in \N}$ such that for each $k \in \N$, and a.e.~$t \in
(0,T)$, $\bm{v}_k(t) \in \andrea{\bY_k}$ and $\bm{v}_k \to \u$
strongly in $L^2(0,T;\X)$.  It follows that the difference $\u_k -
\bm{v}_k$ converges weakly to zero in $L^2(0,T;\X)$ \andrea{as
  $k\to\infty$}.  Moreover, we can consider $\bet = (\u_k - \bm{v}_k)$
in \eqref{G:a:u} and obtain thanks to the coercivity property
\eqref{prop:E:1} that 
\begin{equation}\label{trick}
\begin{aligned}
0 & = (\C (\ph_k)(\E(\u_k)- \bar \E(\ph_k)) , \E(\u_k - \bm{v}_k)) - ( \g, \u_k - \bm{v}_k)_{\GN} \\
& = (W_{,\E}(\varphi_k,\E(\u_k)) - W_{,\E}(\varphi_k,\E(\bm{v}_k)), \E(\u_k - \bm{v}_k)) \\
& \quad + (\C(\ph_k)(\E(\bm{v}_k) - \bar \E(\ph_k)), \E(\u_k - \bm{v}_k)) - (\g, \u_k - \bm{v}_k)_{\GN} \\
& \geq C_\andrea{4} \no{\E(\u_k - \bm{v}_k)}_{L^2}^2 + (\C(\ph_k)(\E(\bm{v}_k) - \bar \E(\ph_k)), \E(\u_k - \bm{v}_k)) - (\g, \u_k - \bm{v}_k)_{\GN}.
\end{aligned}
\end{equation}
Integrating the above inequality over $(0,T)$, and applying the strong
convergence of $\bm{v}_k$ to $\u$, the convergence properties of $\ph_k$
to $\ph$, and the weak convergence of
$\u_k - \bm{v}_k $ to zero in $L^2(0,T;\X)$, we then obtain
\begin{align*}
\no{\E(\u_k - \bm{v}_k)}_{L^2(Q)}^2 \to 0 \quad \text{ as } k \to \infty.
\end{align*}
By Korn's inequality this shows that $\u_k - \bm{v}_k$ converges strongly to zero in $L^2(0,T;\X)$ and hence $\u_k \to \u$ strongly in $L^2(0,T;\X)$ and a.e.~in $Q$.

\andrea{Now, we aim at passing to the limit.}
The standard procedure is to fix $j \in \N$ in \eqref{Gal:a}, multiply \eqref{G:a:phi}, \eqref{G:a:mu}, \eqref{G:a:sig}, and \eqref{G:a:u} with an arbitrary $\theta \in C^{\infty}_c(0,T)$, pass to the limit $k \to \infty$ with the above compactness results, and use the density of $\cup_{k \in \N} Z_k$ in $H^1(\Omega)$, and the density of $\cup_{k \in \N} \bY_k$ in $\X$ to show that the limit $(\varphi,\mu,\sigma,\u)$ satisfies \eqref{w:1}-\eqref{w:4} for all $\zeta \in L^2(0,T;H^1(\Omega))$, $\xi \in L^2(0,T;H^1(\Omega)) \cap L^{\infty}(Q)$ and $\bet \in L^2(0,T;\X)$.  

The first step is \andrea{to} employ the above compactness assertions to pass to the limit $k \to \infty$ and recover \eqref{w:4}.  This can be done since we have the a.e.~convergence and strong convergence of $\varphi_k$ to $\varphi$ and the particular form of $W_{,\E}(c, E)$.  
For the other equations, a key point to 
pass to the limit is the strong convergence of 
$\E(\u_k)$ to $\E( \u)$ in $L^2(0,T;L^2(\Omega)^{d \times d})$ which follows
from the strong convergence of $\u_k$.  We omit the easy details and sketch the less obvious points.

\paragraph{Mobility term.} Continuity of $m$ from \eqref{ass:m} and the a.e.~convergence of $\varphi_k$ (resp.~$\sigma_k$, $\u_k$ and $\E(\u_k)$) to $ \varphi$ (resp.~$ \sigma$, $\u$ and $\E( \u)$) leads to $m(\varphi_k, \sigma_k, \u_k, \E(\u_k)) \to m( \varphi,  \sigma, \u, \E( \u))$ a.e.~in $Q$.  Boundedness of $m$ from \eqref{ass:m} and the dominated convergence theorem yields that 
\begin{align*}
\theta(t) m(\varphi_k, \sigma_k, \u_k,  \E(\u_k)) \nabla z_j \to \theta(t) m( \varphi,  \sigma, \u, \E( \u)) \nabla z_j \quad \mbox{ strongly in } L^2(Q),
\end{align*}
so that together with the weak convergence of $\nx \mu_k$ to $\nx  \mu$ in $L^2(Q)$ we find that 
\begin{align*}
\int_0^T \theta(t) (m(\varphi_k, \sigma_k, \u_k,  \E(\u_k)) \nx \mu_k, \nx z_j) \dt \to \int_0^T \theta(t) (m( \varphi,  \sigma, \u,  \E( \u)) \nx  \mu, \nx z_j) \dt.
\end{align*}
A similar argument can be used to show the convergence involving
$\tilde U(\varphi_k, \sigma_k, \E(\u_k))$ using the boundedness and
continuity of $\tilde U$ with respect to its arguments.

\paragraph{Potential term.}  Continuity of $\psi'$ and the a.e.~convergence of $\varphi_k$ to $ \varphi$ in $Q$ yields $\psi'(\varphi_k) \to \psi'( \varphi)$ a.e.~in $Q$.  The above compactness results, the sublinear growth of $\psi_2'$ and the generalised dominated convergence theorem lead to
\begin{align*}
\psi_2'(\varphi_k) \to \psi_2'( \varphi) \quad \text{ strongly in } L^2(Q).
\end{align*}
For the monotone part $\psi_1'(\varphi_k)$, we show that the family
$\{\theta(t) \psi_1'(\varphi_k) z_j\}_{k \in \N}$ is uniformly
integrable over $Q$, so that together with the a.e.~convergence
$\psi_1'(\varphi_k)$ to $\psi_1'( \varphi)$ in $Q$, we obtain via
Vitali's convergence theorem that
\begin{align*}
\int_{Q} \theta(t) \psi_1'(\varphi_k) z_j d\x \dt \to \int_{Q} \theta(t) \psi_1'( \varphi) z_j d\x \dt.
\end{align*}
We now show the uniform integrability.
Let $\eta > 0$ be arbitrary, then choosing $\delta > 0$ sufficiently small so that\begin{align*}
\delta \Big ( T C(1+\beta) + C_{\delta} \Big ) \no{\theta}_{L^{\infty}(0,T)} \no{z_j}_{H^2(\Omega)} < \eta,
\end{align*}
where $C$ is the constant in \eqref{Est1}, and $C_\delta$ is the constant in \eqref{ass:pot} associated to $\delta$, we obtain from \eqref{ass:pot} and the fact $z_j \in H^2(\Omega)$ that for any measurable subset $E \subset Q$ with $\abs{E} < \delta$,
\begin{align}
\notag \int_{E} \abs{ \theta(t) \psi_1'(\varphi_k) z_j } d\x \dt & \leq \no{\theta}_{L^{\infty}(0,T)} \no{z_j}_{H^\andreaold{2}(\Omega)} \int_{E} \Big ( \delta \psi_1(\varphi_k) + C_{\delta} \Big ) d\x \dt \\
\label{psicon}
& \leq \Big (\delta \andrea{T} C(1+\beta)  + C_{\delta} \abs{E}\Big ) \no{\theta}_{L^{\infty}(0,T)} \no{z_j}_{H^2(\Omega)} < \eta,
\end{align}
which implies the uniform integrability of the family $\{ \theta(t) \psi_1'(\varphi_k) z_j\}_{k \in \N}$.

\paragraph{Elasticity terms.}
Thanks to the strong convergence of $\u_k$ in $L^2(0,T;\X)$
and of $\ph_k $ in $\L2 {\LOx2}$
\andrea{and due to \eqref{prop:E:2} and \eqref{prop:E:3}}, we have
\begin{align*}
C_\andrea{5} (1 + \abs{\varphi_k}^2 + \abs{\E(\u_k)}^2) &\to C_\andrea{5} (1 + \abs{ \varphi}^2 + \abs{\E( \u)}^2) 
\quad \text{strongly in } L^1(Q), \\
C_\andrea{5} (1 + \abs{\varphi_k} + \abs{\E(\u_k)}) &\to C_\andrea{5} (1 + \abs{ \varphi} + \abs{\E( \u)}) \qquad  \text{strongly in } L^2(Q),
\end{align*}
and so by \eqref{ass:W} and the generalised dominated convergence theorem, we infer that
\begin{align*}
\int_Q \theta(t) W_{,\E}(\varphi_k, \E(\u_k)) \cdot \nabla \y_j d\x \dt & 
\to \int_{Q} \theta(t) W_{,\E}( \varphi, \E( \u)) \cdot \nabla \y_j d\x \dt, \\
\int_Q \theta(t) W_{,\varphi}(\varphi_{k}, \E(\u_k)) z_j d\x \dt &
\to \int_Q \theta(t) W_{,\varphi}( \varphi, \E( \u)) z_j d\x \dt
\end{align*}
on account of the fact that $z_j \in H^2_{\mathbf{n}}(\Omega) \subset L^{\infty}(\Omega)$.

\paragraph{Comparison principle.} To establish the boundedness property \eqref{sigma:bdd} for $ \sigma$, so that $\tilde U( \varphi,  \sigma, \E( \u)) = U( \varphi,  \sigma, \E( \u))$, i.e., $g( \sigma) =  \sigma$, we employ a comparison principle.  We recall that the positive part $f_{+}$ and negative part $f_{-}$ of a function $f$ are defined as
\begin{align*}
f_{+}(x) = \max(f(x), 0), \quad f_{-}(x) = \max(-f(x), 0),
\end{align*}
so that $f(x) = f_{+}(x) - f_{-}(x)$.  Testing \eqref{w:3} with $- (
\sigma)_{-}$ and using the relations
\begin{align*}
\inn{ \sigma_t}{( \sigma)_{-}} = - \frac{d}{dt} \no{( \sigma)_{-}}_{L^2}^2, \quad ( \nx  \sigma, \nx ( \sigma)_{-}) = - \no{\nx ( \sigma)_{-}}_{L^2}^2, 
\end{align*}
we obtain
\begin{align*}
& \frac{\beta}{2} \frac{d}{dt} \no{( \sigma)_{-}}_{L^2}^2 + \no{\nx ( \sigma)_{-}}_{L^2}^2 + \kappa \no{( \sigma)_{-}}_{L^2_\Gamma}^2 + \lambda_c \no{h^{1/2}( \sigma)_{-}}_{L^2}^2 + B \no{( \sigma)_{-}}_{L^2}^2  \\
& \quad = - B(\sigma_c, ( \sigma)_{-}) - \kappa (\sigma_B, ( \sigma_{-}))_\Gamma \leq 0
\end{align*}
on account of the fact that $\sigma_B, \sigma_c$ and $( \sigma)_{-}$
all are non-negative.  Integrating the above inequality and using the fact
that $\sigma_0$ is non-negative we obtain
\begin{align*}
\no{( \sigma)_{-}(t)}_{L^2}^2 \leq \no{(\sigma_0)_{-}}_{L^2}^2 \leq 0 \quad \text{ for all } t \in (0,T),
\end{align*}
so \andrea{that} $ \sigma$ is non-negative a.e.~in $Q$.  On the other
hand, testing the equation \eqref{w:3} by $( \sigma - M)_{+}$, where
we recall that $M = \max( \no{\sigma_c}_{L^{\infty}(Q)},
\no{\sigma_B}_{L^{\infty}({\Sigma})})$,  yields
\begin{align*}
& \frac{\beta}{2} \frac{d}{dt} \no{( \sigma - M)_{+}}_{L^2}^2 + \no{\nx ( \sigma - M)_{+}}_{L^2}^2 + \kappa \no{( \sigma - M)_{+}}_{L^2_\Gamma}^2 \\
& \qquad + \lambda_c \no{h^{1/2}( \sigma - M)_{+}}_{L^2}^2 + B \no{( \sigma - M)_{+}}_{L^2}^2 \\
& \quad = - \kappa (M - \sigma_B, ( \sigma - M)_{+})_{\Gamma} - \lambda_c (h( \varphi) M, ( \sigma - M)_{+}) - B(M - \sigma_c, ( \sigma - M)_{+}) \\
& \quad \leq 0
\end{align*}
on account of the fact that $h$, $M - \sigma_c$, $M - \sigma_B$ and $(
\sigma - M)_{+}$ all are non-negative.  Integrating the above inequality and using the fact that $\sigma_0 \leq M$ a.e.~in $\Omega$ leads to
\begin{align*}
\no{( \sigma - M)_{+}(t)}_{L^2}^2 \leq \no{(\sigma_0 - M)_{+}}_{L^2}^2 = 0 \quad \text{ for all } t \in (0,T)
\end{align*}
\andrea{so that} $  \sigma \leq M$ a.e.~in $Q$ as we claimed.

\subsection{Energy inequality}
Let us now combine inequalities \eqref{Est1}, \eqref{Est2} and
\eqref{Est3} to obtain that
\begin{equation*}
\begin{aligned}
& \sup_{t \in (0,T)} \Big ( \no{\varphi_k(t)}_{H^1}^2 
+ \no{\psi(\varphi_k(t))}_{L^1} + \beta \no{\sigma_k(t)}_{L^2}^2 + \no{\u_k(t)}_{\X}^2 \Big ) \\
& \qquad + \no{\mu_k}_{L^2(0,T;H^1)}^2 
+ \no{\sigma_k}_{L^2(0,T;H^1)}^2  
+ \no{\ph_k}_{H^1(0,T;(H^1)')}^2  
+ \beta\no{\sigma_k}_{H^1(0,T;(H^1)')}^2  
\\ & \quad \leq
C\Big ( 1 + \no{\psi(\varphi_k)}_{L^1(0,t;L^1)} + \no{\varphi_k}_{L^2(0,t;H^1)}^2 + \no{\E(\u_k)}_{L^2(0,t;L^2)}^2 \Big ) \\
& \qquad + C \no{\varphi_{k,0}}_{H^1}^2 
+ \no{\psi(\varphi_{k,0})}_{L^1} + \beta \no{\sigma_{k,0}}_{L^2}^2.
\end{aligned}
\end{equation*}
%  
%  
%We return to the inequality \eqref{Energy:ptl}, and set $K = 0$ (so that the last term on the right-hand side involving $\no{\nabla \varphi_k}_{L^2}^2$ drops out), integrating in time from $0$ to $s \in (0,T)$, and using \andrea{\eqref{Poin} and }\eqref{W:0:est}, leads to
%\begin{equation}\label{Energy:ptl:1}
%\begin{aligned}
%& \no{\nabla \varphi_k(s)}_{L^2}^2 + \no{\psi(\varphi_k(s))}_{L^1} + \beta \no{\sigma_k(s)}_{L^2}^2 
%\\
%& \qquad + \int_\Omega W(\varphi_k(s), \E(\u_k(s))d\x - \int_{\GN} \g \cdot \u_k(s)\dA 
%+ \int_0^s 
%\no{\nabla \mu_k}_{L^2}^2 + \no{\sigma_k}_{H^1}^2 \dt \\
%& \quad \leq  C \int_0^s \abs{\mean{\mu_k}} - \no{\sigma_k}_{L^2}^2 + \no{\varphi_k}_{L^2}^2 \dt \\
%& \qquad + C \Big (1 +  \no{\varphi_{k,0}}_{H^1}^2  \andrea{+} \no{\psi(\varphi_{0,k})}_{L^1} + \beta \no{\sigma_{0,k}}_{L^2}^2  \Big ).
%\end{aligned}
%\end{equation}
From the compactness assertions stated in Section \ref{sec:compact},
we infer, by Fatou's lemma and the non-negativity of $\psi$ 
that, for $a.e. s \in (0,T)$, it holds 
\begin{align*}
	\int_\Omega \psi(\ph(s)) \dx 
	\leq \liminf_{k \to \infty } \int_\Omega \psi(\ph_k(s)) \dx. 
\end{align*}
Moreover, invoking the weak/weak* lower semicontinuity of the norms, using the properties
$\no{\varphi_{k,0}}_{H^1}^2 \leq C \no{\varphi_{0}}_{H^1}^2$,
$\no{\sigma_{k,0}}_{L^2}^2 \leq \no{\sigma_0}_{L^2}^2$ originating
from the orthogonality {in $\Hx1$} of the basis functions
$\{z_i\}_{i \in \N}$, and recalling \eqref{psi:ini:conv},
allow us to pass to the limit as $k\to\infty$ in the above inequality 
to obtain \eqref{En:ineq}.

\subsection{More general initial conditions}\label{sec:ini}
To complete the proof of Theorem \ref{thm:exist}, we now assume that
$\varphi_{0} \in H^1(\Omega)$ with $\psi(\varphi_0) \in L^1(\Omega)$
and use ideas of \cite{CFG}.
For any $\delta \in (0,1]$, we denote by
$\varphi_{0, \delta} \in H^2_{\mathbf{n}}(\Omega)$ the unique solution to the
elliptic problem:
\begin{align}\label{Ell:ini}
\begin{cases}
- \delta \Lx \varphi_{0,\delta} + \varphi_{0,\delta} = \varphi_{0} & \text{ in } \Omega, \\
\pdnu \varphi_{0,\delta} = 0 & \text{ on } \Gamma.
\end{cases}
\end{align}
The well-posedness and regularity of $\varphi_{0,\delta}$ follows from standard application of the Lax--Milgram theorem and elliptic regularity theory 
{(see, e.g., \cite{Gris} for 
the corresponding regularity theory for convex domains).}
Furthermore, testing the above equation by $\varphi_{0,\delta}$ and $- \Lx \varphi_{0,\delta}$, respectively leads to the uniform estimates
\begin{equation}\label{varphidelta:H1}
\begin{aligned}
2\delta \no{\nx \varphi_{0,\delta}}_{L^2}^2 + \no{\varphi_{0,\delta}}_{L^2}^2 & \leq \no{\varphi_{0}}_{L^2}^2, \\
2 \delta \no{\Lx \varphi_{0,\delta}}_{L^2}^2 + \no{\nx \varphi_{0,\delta}}_{L^2}^2 & \leq \no{\nabla \varphi_{0}}_{L^2}^2.
\end{aligned}
\end{equation}
This implies that, as $\delta$ goes to zero, we have
\begin{align*}
\varphi_{0,\delta}&  \to \varphi_0 \quad \mbox{ weakly in } H^1(\Omega), \\
\varphi_{0,\delta}& \to \varphi_0 \quad \mbox{ strongly in } L^2(\Omega) \text{ and a.e. in } \Omega.
\end{align*}
From \eqref{ass:pot}, we see that the function 
\begin{align*}
G(s) := \psi(s) + \tfrac{1}{2} C_1 \abs{s}^2
\end{align*}
is convex and non-negative, since
\begin{align*}
G''(s) = \psi_{1}''(s) + \psi_{2}''(s) + C_1 \geq 0.
\end{align*}
Furthermore, by the assumption \eqref{ass:ini} on $\varphi_0$ it holds
that $G(\varphi_0) \in L^1(\Omega)$.  Then, testing the elliptic
problem \eqref{Ell:ini} with $G'(\varphi_{0,\delta})$ yields
\begin{align*}
\int_{\Omega} (\varphi_{0,\delta} - \varphi_0) G'(\varphi_{0,\delta}) d\x = - \int_{\Omega} \delta G''(\varphi_{0,\delta}) \abs{\nx \varphi_{0,\delta}}^2 d\x \leq 0.
\end{align*}
Since $\psi'(0) = 0$, we see that $G'(0) = 0$, and using the convexity
of $G$ and the previous inequality we infer that
\begin{align*}
\int_{\Omega} G(\varphi_{0,\delta}) d\x \leq \int_{\Omega} G(\varphi_0) + G'(\varphi_{0,\delta})(\varphi_{0,\delta} - \varphi_0) d\x \leq \int_{\Omega} G(\varphi_0) d\x < \infty.
\end{align*}
In particular, by the strong convergence of $\varphi_{0,\delta} \to \varphi_0$ in $L^2(\Omega)$ and the weak lower semicontinuity of the $L^2(\Omega)$-norm, it holds that 
\begin{equation}\label{psi:delta:L1}
\begin{aligned}
\limsup_{\delta \to 0} \int_{\Omega} \psi(\varphi_{0,\delta}) d\x &
\leq \limsup_{\delta \to 0} \int_{\Omega} G(\varphi_{0,\delta}) d\x +
\limsup_{\delta \to 0} \int_\Omega - \frac{C_1}{2} \abs{\varphi_{0,\delta}}^2 d\x\\
&  \leq \int_{\Omega} G(\varphi_0)\andrea{d\x} - \lim_{\delta \to 0} \int_\Omega \frac{C_1}{2} \abs{\varphi_{0,\delta}}^2 d\x \leq \int_\Omega \psi(\varphi_0) d\x.
\end{aligned}
\end{equation}
Hence, for given $\varphi_0 \in H^1(\Omega)$ satisfying
$\psi(\varphi_0) \in L^1(\Omega)$, we consider the sequence of
solutions $(\varphi_\delta, \mu_\delta, \sigma_\delta, \u_\delta)$ to
\eqref{Lima} with initial conditions $(\varphi_{0,\delta}, \sigma_0)$
such that $\varphi_{0,\delta} \in H^2_{\mathbf{n}}(\Omega)$ is the
  unique solution to \eqref{Ell:ini}. 
  Then,
  $(\varphi_\delta, \mu_\delta, \sigma_\delta, \u_\delta)$ satisfies
  the energy inequality \eqref{En:ineq} with a right-hand side given by
\begin{align*}
 &  C\Big ( 1 + \no{\psi(\varphi_{\delta})}_{L^1(0,t;L^1)} + \no{\varphi_\delta}_{L^2(0,t;H^1)}^2 + \no{\E(\u_\delta)}_{L^2(0,t;L^2)}^2 \Big ) \\
& \quad + C \no{\varphi_{\delta,0}}_{H^1}^2 
+ \no{\psi(\varphi_{\delta,0})}_{L^1} + \beta \no{\sigma_{0}}_{L^2}^2
\end{align*}
with a positive constant $C$ independent of $\delta \in (0,1]$.
Then, a Gronwall argument yields the uniform in $\delta$ estimate so that the solution $(\varphi_\delta, \mu_\delta, \sigma_\delta, \u_\delta)$ satisfies the same compactness assertions listed in Section \ref{sec:compact} and converges along a non-relabelled subsequence to limit functions $(\varphi, \mu, \sigma, \u)$ in the limit $\delta \to 0$.  The strong convergence of $\E(\u_{\delta}) \to \E(\u)$ follows from the monotonicity argument outlined in \cite{GHab,Gela}.  We omit the rest of the details and infer that $(\varphi, \mu, \sigma, \u)$ is a weak solution fulfilling the assertions of Theorem \ref{thm:exist} with initial condition $\varphi_0$ satisfying \eqref{ass:ini}.

\section{Quasi-static limit}
We now consider the limit $\beta\to 0$ in
\eqref{Lima}, i.e. we consider the quasi-static limit of the
nutrient diffusion equation. In order to prove
Theorem~\ref{thm:quasi}, we denote the solutions to the
weak formulation \eqref{w:1}-\eqref{w:4} constructed in
Theorem~\ref{thm:exist} by $(\varphi_\beta, \mu_\beta, \sigma_\beta,
\u_\beta)$. 
The compactness assertions aside from \eqref{sigma:beta} are
consequences of the uniform estimates obtained from analogues of
\eqref{En:ineq}, as well as the monotonicity argument of
\cite{GHab,Gela}.  From the uniform estimates we also have
$\sigma_\beta \to \sigma_*$ weakly in $L^2(0,T;H^1(\Omega))$, which is
sufficient to pass to the limit in the source terms
$U(\varphi_\beta, \sigma_\beta, \E(\u_\beta))$ and
$S(\varphi_\beta, \sigma_\beta)$ due to their particular forms
\eqref{defn:U} and \eqref{defn:S}, as well as the strong convergence
of $\varphi_\beta \to \varphi_*$ in $C^0([0,T]; L^r(\Omega))$
\andrea{for any $r<\infty$ if $d=2$ or $r<6$ if $d=3$.} However, since
the mobility $m$ depends (perhaps non-linearly) on $\sigma$, we
require the a.e.~convergence of $\sigma_\beta$ to $\sigma_*$ in $Q$
which is not available simply from the uniform $L^2(0,T;H^1(\Omega))$
estimate for $\{\sigma_\beta\}_{\beta \in (0,1]}$.  Therefore, in the
following we derive a strong convergence result for $\sigma_\beta$ in
$L^2(0,T;L^2(\Omega))$.  First, considering $\zeta\in
L^2(0,T;H^1(\Omega))$, and passing to the limit $\beta \to 0$ in
\begin{align*}
0 = \int_0^T \beta \inn{\sigma_{\beta,t}}{\zeta} + (\nabla \sigma_{\beta},\nabla \zeta) + \kappa (\sigma_\beta - \sigma_B, \zeta)_{\Gamma} - (S(\varphi_\beta, \sigma_\beta), \zeta) \dt
\end{align*}
yields
\begin{align*}
0 = \int_0^T (\nabla \sigma_*,\nabla \zeta) + \kappa (\sigma_* - \sigma_B, \zeta)_{\Gamma} - (S(\varphi_*, \sigma_*), \zeta) \dt.
\end{align*}
Then, denoting $\hat \sigma:= \sigma_\beta - \sigma_*$ and taking the difference of the two equations above gives
\begin{equation}\label{str:sigma:beta}
\begin{aligned}
& \int_0^T (\nabla \hat \sigma, \nabla \zeta) + \kappa (\hat \sigma, \zeta)_{\Gamma} + B(\hat \sigma, \zeta) + \lambda_c (h(\varphi_*) \hat \sigma, \zeta) \dt \\
& \quad = - \int_0^T \inn{\beta \sigma_{\beta, t}}{\zeta} - \lambda_c ([h(\varphi_{\beta}) - h(\varphi_*)] \sigma_\beta, \zeta) \dt.
\end{aligned}
\end{equation}
Choosing $\zeta= \sigma_\beta - \sigma_*$ and observe that
\begin{align*}
\abs{\int_0^T \inn{\beta \sigma_{\beta,t}}{\sigma_\beta - \sigma_*} \dt} & = \abs{\frac{\beta}{2} \no{\sigma_{\beta}(T)}_{L^2}^2 - \frac{\beta}{2} \no{\sigma_0}_{L^2}^2 - \int_0^T \inn{\beta \sigma_{\beta, t}}{\sigma_*} \dt} \\
& \leq \frac{\beta}{2} M + \frac{\beta}{2} \no{\sigma_0}_{L^2}^2 + \abs{\int_0^T \inn{\beta \sigma_{\beta, t}}{\sigma_*} \dt} \to 0, 
\end{align*}
and
\begin{align*}
\int_0^T \lambda_c ([h(\varphi_\beta) - h(\varphi_*)] \sigma_\beta, \sigma_\beta - \sigma_*) \dt & \leq \max_{t \in [0,T]} \lambda_c(t) C \no{h(\varphi_{\beta}) - h(\varphi_*)}_{L^2(0,T;L^2)} \to 0,
\end{align*}
on account of the weak convergence of $\beta \sigma_{\beta, t} \to 0$
in $L^2(0,T;H^1(\Omega)')$, the strong convergence
$h(\varphi_{\beta}) \to h(\varphi_*)$ in $L^2(0,T;L^2(\Omega))$ (which
is a consequence of the boundedness of $h$ and the a.e.~convergence of
$\varphi_\beta \to \varphi$ in $Q$), as well as the boundedness
$0 \leq \sigma_{\beta}, \sigma_* \leq M$ a.e.~in $Q$.

Hence, choosing $\zeta = \sigma_\beta - \sigma_*$ in \eqref{str:sigma:beta}, neglecting the non-negative term $B \no{\hat \sigma}_{L^2}^2 + \lambda_c \no{h^{1/2}(\varphi_*) \hat \sigma}_{L^2}^2$ on the left-hand side, and employing the generalised Poincar\'{e} inequality \eqref{Poin} yields
\begin{align*}
\no{\sigma_\beta - \sigma_*}_{L^2(0,T;H^1)} \to 0 \mbox{ as } \quad \beta \to 0.
\end{align*}
This yields the compactness assertion \eqref{sigma:beta}, and the rest of the proof follows similarly as described in the previous sections.

\section{Regularity}

\subsection{Regularity for the nutrient}
Suppose $\sigma_0 \in H^1(\Omega)$ and $\sigma_B \in H^1(0,T;L^2(\Gamma))$.  The following formal estimates can be obtained rigorously at the level of the Galerkin approximation, and so we will only sketch the details.  The nutrient system can be expressed as
 \begin{align*}
\begin{cases}
\beta \sigma_t - \Lx \sigma + B\s = - \lambda_c h(\ph) \s + B \s_c =:\andrea{f_\s} & \text{ in } Q, \\
\pdnu \sigma + \kappa (\sigma - \sigma_B)  = 0	& \text{ on } \Sigma,\\
\s(0)=\s_0	& \text{ in }\Omega.
\end{cases}
\end{align*}
From Theorem~\ref{thm:exist} and from the assumption on the data, it easily follows that $\andrea{f_\s}\in L^2(Q)$.  Hence, testing by $\s_t$ yields
\begin{align*}
%\frac{d}{dt} \frac{1}{2} \Big ( \beta \no{\nabla \s}_{L^2}^2 + \beta \kappa \no{\s}_\andrea{{L^2_\Gamma}}^2 \Big ) + \no{\Delta \s}_{L^2}^2 &\leq \no{\andrea{f_\s} - B \s}_{L^2} \no{\Delta \s}_{L^2}
%\andreaold{+ \kappa\int_\Gamma \beta  \s_B \s_t \dA}
%, \\
  \frac{d}{dt} \frac{1}{2} \andrea{\Big(} \no{\nabla \s}_{L^2}^2 + B \no{\s}_{L^2}^2 + \kappa \no{\s}_{\andrea{{L^2_\Gamma}}}^2 \Big ) + \beta \no{\s_t}_{L^2}^2 &\leq \no{\andrea{f_\s}}_{L^2} \no{\s_t}_{L^2} + \andrea{\kappa}\int_\Gamma  \s_B \s_t \dA.
\end{align*}

In order to handle the last boundary term 
we integrate in time and by parts, and invoke the Young inequality and
the trace theorem to obtain that
\begin{align*}
	& \andrea{\kappa}\int_0^t \int_\Gamma \s_B \s_t \dA
	\\  & \quad= - \andrea{\kappa}\int_0^t \int_\Gamma \s_{B,t} \, \s \dA
	+\andrea{\kappa}\int_\Gamma \s_B(t) \s(t) \dA
	-\andrea{\kappa}\int_\Gamma \s_B(0) \s_0 \dA
	\\ & \quad \leq 
	\andrea{\frac \d 4} \int_\Gamma |\s(t)|^2
	+C (\no{\s_{B,t}}_{L^2(\Sigma)}^2 
	+ \no{\s}_{L^2(\Sigma)}^2 
	%+\no{\s_B \s}_{L^\infty(0,T;L^1(\Gamma))}
	+ \no{\s_{B}}_{C^0(0,T;L^2(\Gamma))}^2 
	+ \no{\s_0}_{H^1}^2),
\end{align*}
\andrea{for a positive $\d$ to be chosen as 
\begin{align*}
	\d = 
	\begin{cases}
	 {\kappa}\quad &\hbox{if} \quad \kappa > 0,
	\\
	 {\varrho B}	\quad &\hbox{if} \quad \kappa=0, 
	\end{cases}
\end{align*}
where $\varrho$ denotes the inverse of the trace constant and} we also owe to the improved smoothness of the data 
$\sigma_B \in \H1 {L^2(\Gamma)) \subset C^0(0,T;L^2(\Gamma)}$.
\andrea{Thus, integrating in time shows that 
\begin{align*}
	\sigma \in L^\infty(0,T;H^1(\Omega)) \quad \forall\, \beta \in [0,\infty), 
	\quad \sigma \in \H1 {L^2(\Omega)} \quad \forall\, \beta >0.
\end{align*}
Hence we can now absorb the possible term $\beta \s_t$ in the source
contribution $f_\s$ which will still belong to $L^2(Q)$.  Then,
provided we require $\s_B \in L^2(0,T;H^{1/2}(\Gamma))$, we can invoke
elliptic regularity theory to conclude, independently of $\beta$, that
$\sigma \in L^2(0,T;H^2(\Omega))$ (see, e.g. \cite[Thm.~4.18,
pp.~137-138]{McL}).}

%For the case $\beta = 0$, testing the nutrient equation with $- \Delta \s$ and employing elliptic regularity yields $\s \in L^2(0,T;H^2(\Omega))$.  Meanwhile, testing with $\s$ gives
%\begin{align*}
%\no{\nabla \s}_{L^2}^2 + \kappa \no{\s}_\andrea{{L^2_\Gamma}}^2 \leq \no{\andrea{f_\s}}_{L^2} \no{\s}_{L^2} + \kappa \no{\s_B}_\andrea{{L^2_\Gamma}} \no{\s}_\andrea{{L^2_\Gamma}}.
%\end{align*}
%By the generalised Poincar\'e inequality, as well as the fact that $\andrea{f_\s} \in L^{\infty}(0,T;L^2(\Omega))$ and $\s_B \in L^\infty(0,T;L^2(\Gamma))$, we infer that $\s \in L^\infty(0,T;H^1(\Omega))$.

\subsection{Regularity under Vegard's law and homogeneous elasticity}
Under Vegard's law \eqref{Vegard} and homogeneous elasticity \eqref{C:cons}, the partial derivative $W_{,\ph}(\varphi, \E(\u))$ assumes the following form
\begin{align*}
W_{,\ph}(\ph, \E(\u)) = - \C (\E(\u) - \hat \E - \E^* \ph ) : \E^*,
\end{align*}
and by the regularities stated in Theorem \ref{thm:exist}, we see that $W_{,\ph}(\ph ,\E(\u)) \in L^{\infty}(0,T;L^2(\Omega)^{d \times d})$.  Hence,  it holds that $f := \chi \sigma + \mu - W_{,\ph}(\ph, \E(\u)) - \eps^{-1} \psi_2'(\ph)$ belongs to $L^2(Q)$.  For $N \in \N$, we introduce the truncation
\begin{align*}
\ph_N := \max(-N, \min(N, \ph)) \text{ a.e.~in } Q.
\end{align*}
Then, it is clear that $\ph_N \to \ph$ a.e.~in $Q$, and as $\ph_N$ is bounded it holds $\psi_1'(\ph_N) \in L^2(0,T;H^1(\Omega)) \cap L^{\infty}(Q)$.  Testing \eqref{w:2} with $\psi_1'(\ph_N)$ and using the convexity of $\psi_1$ leads to
\begin{align*}
\int_\Omega \eps^{-1} \psi_1'(\ph) \psi_1'(\ph_N) d\x \leq \int_\Omega \eps^{-1} \psi_1'(\ph) \psi_1'(\ph_N) + \eps \psi_1''(\ph_N) \abs{\nabla \ph_N}^2 d\x \leq \no{f}_{L^2} \no{\psi_1'(\ph_N)}_{L^2}.
\end{align*} 
Using the facts that $\psi_1'$ is increasing and the monotonicity
$\psi_1'(s) s\geq 0$, we infer
\begin{align*}
\no{\psi_1'(\ph_N)}_{L^2}^2 \leq \int_\Omega \psi_1'(\ph) \psi_1'(\ph_N) d\x \leq \eps \no{f}_{L^2} \no{\psi_1'(\ph_N)}_{L^2}.
\end{align*}
This gives boundedness of $\psi_1'(\ph_N)$ in $L^2(Q)$, and by Fatou's lemma,
\begin{align*}
\no{\psi_1'(\ph)}_{L^2(Q)} \leq \liminf_{N \to \infty} \no{\psi_1'(\ph_N)}_{L^2(Q)} \leq \eps \no{f}_{L^2(Q)},
\end{align*}
which is the first assertion.  In turn, \eqref{w:2} is the weak formulation of the elliptic problem
\begin{align*}
\begin{cases}
- \eps \Lx \ph = f - \eps^{-1} \psi_1'(\ph) & \text{ in } \Omega, \\
\pdnu \ph = 0 & \text{ on } \Gamma,
\end{cases}
\end{align*}
and elliptic regularity then yields
\begin{align*}
\int_0^T \no{\ph}_{H^2}^2 \dt \leq \int_0^T  C \Big ( \no{\ph}_{H^1}^2 + \no{f}_{L^2}^2 + \eps^{-1} \no{\psi_1'(\ph)}_{L^2}^2 \Big ) \dt < \infty
\end{align*}
which completes the proof.

\andreaold{\begin{remark}
Having $\ph\in \L2 {H^2}$ at disposal, it would be natural to hope
to improve the regularity of the displacement $\u$ since, roughly speaking,
equation \eqref{u} can be written as an elliptic equation for $\u$ 
whose source has a regularity which depends on $\ph$. 
Unfortunately, this is not in general possible due to
the choice of the boundary conditions which prevent the regularity of $\u$
to be improved.
\end{remark}}

\section{Continuous dependence under Vegard's law and homogeneous elasticity}
%This section is devoted to show some continuous dependence results, which
%as a consequence imply the uniqueness of the solution.
%Before moving to the result, let us present some tools that will be important in this section.
\subsection{Preliminaries}
The Riesz isomorphism $\A: H^1(\Omega) \to H^1(\Omega)'$ is defined as
\begin{align*}
\inn{\A u}{v} := (u,v)_{H^1} = \int_\Omega (\nabla u \cdot \nabla v + u v) d\x.
\end{align*}
It is well-known that the restriction of $\A$ to $H^2_{\mathbf{n}}(\Omega)$
\andrea{yields} an isomorphism from $H^2_{\mathbf{n}}(\Omega)$ to
$\LOx2$, so that its inverse
$\A^{-1} : \LOx2 \to H^2_{\mathbf{n}}(\Omega)$ is well-defined.
Moreover, the following properties hold
\begin{align*}
\inn{\A u}{\A^{-1}w} &= \inn{w}{u} \quad \text{ for all } u \in H^1(\Omega), \; w \in H^1(\Omega)', \\
\inn{w}{\A^{-1}z} &= (w, z )_{*} \quad \text{ for all } w,z \in H^1(\Omega)',
\end{align*}
where the symbol $(\cdot, \cdot )_{*}$ denotes the standard 
inner product in the dual of $H^1(\Omega)$, and $\inn{\cdot}{\cdot}$ is the duality pairing between $H^1(\Omega)$ and its dual.  By Lax--Milgram theorem, we have the estimate
\begin{align*}
\no{\A^{-1} z}_{H^1}^2 = \inn{z}{\A^{-1} z} \leq \no{z}_* \no{\A^{-1} z}_{H^1} \quad \implies \quad \no{\A^{-1} z}_{H^1} \leq \no{z}_{*}
\end{align*}
for all $z \in H^1(\Omega)'$.  Furthermore, the continuous embedding  $H^1(\Omega) \subset L^2(\Omega)$ yields
\begin{align}
\inn{u}{v} &= \int_\Omega u v \quad \text{ for all } u,v\in L^2(\Omega), \label{identification}
\end{align}
which implies, for any $f \in H^1(\Omega)$, the following interpolation inequality
\begin{align}\label{interpol_one}
\no{f}_{L^2}^2 = \inn{f}{f} = \inn{\A f}{\A^{-1}f} \leq \no{f}_{H^1} \no{\A^{-1} f}_{H^1} \leq \no{f}_{H^1} \no{f}_*.
\end{align}
Moreover, for all $w \in \H1 {H^1(\Omega)'}$, we also infer that
\begin{align*}
\inn{w_t(t)}{\A^{-1}w(t)} = \frac 12 \frac d {dt} \norma{w(t)}^2_{*} 
\quad \text{ for a.e~} t \in (0,T).
\end{align*}
We also state here two special cases of the Gagliardo--Nirenberg inequality 
in three dimensions that we will use later:
\begin{align}
\no{f}_{L^3} &\leq  C \no{f}_{L^2}^{1/2} \no{f}_{H^1} ^{1/2} \quad \text{ for all } f \in H^1(\Omega), \label{interpol_two}\\
\no{f}_{L^\infty} &\leq C \no{f}_{H^1}^{1/2}\no{f}^{1/2}_{H^2} \quad \text{ for all } f \in H^2(\Omega). \label{GN:3d}
\end{align}
We will prove the continuous dependence result for $d=3$.
The case $d=2$ is easier due to better embedding properties and is omitted.

Taking the difference of system \eqref{Lima} between two sets of solutions $\{(\ph_i, \mu_i, \s_i, \u_i)\}_{i=1,2}$, and denoting the differences as
\begin{align}
\ph:= \ph_1-\ph_2, \quad \mu:= \mu_1-\mu_2, \quad \s:= \s_1-\s_2, \quad \u:= \u_1-\u_2,\label{diff_one}
\end{align}
along with 
\begin{align}
\notag f_i &:= f(\varphi_i), \quad h_i := h(\varphi_i), \quad k_i := k(\varphi_i), \quad \psi_i' := \psi'(\ph_i), \quad W_{,\E,i} := W_{,\E}(\ph_i, \E(\u_i)), \\\notag 
\hat f &:= f_1 - f_2, \quad 	\hat h := h_1 - h_2, \quad 	\hat k := k_1 - k_2, \quad \hat \psi' := \psi_1' - \psi_2' , \\\notag 
\hat \sigma_B &:= \sigma_{B,1} - \sigma_{B,2}, \quad \hat \sigma_c := \sigma_{c,1} - \sigma_{c,2}, \quad \hat \g := \g_1 - \g_2, \\
\hat{W}_{,\ph} & := W_{,\ph}(\ph_1, \E(\u_1)) - W_{,\ph}(\ph_2, \E(\u_2)) = - \C(\E(\andrea{ \u}) - \E^* \ph): \E^*,
\label{diff_two}
\end{align}
for $i = 1,2$,  it holds that
\begin{subequations}
\begin{alignat}{2} 
\notag 0 & = \inn{\ph_t}{\zeta} + (\nabla \mu,  \nabla  \zeta) - \andrea{\Big(}\tfrac{\lambda_p}{1+|W_{,\E,1}|} ( \hat f \s_1 + f_2 \s)  - \lambda_a \hat k, \zeta\andrea{\Big)} \\
\label{cts:1}  & \quad - \Big (\tfrac{\lambda_p f_2 \s_2}{(1+|W_{,\E,1}|)(1+|W_{,\E,2}|)} (|W_{,\E,1}|-|W_{,\E,2}|), \zeta \Big ), \\ 
\label{cts:2} 0 & = (\mu, \xi) - \eps (\nabla \ph,\nabla   \xi)  - \eps^{-1} (\hat \psi', \xi) + \chi (\sigma, \xi) + (\C (\E(\u) - \E^* \ph): \E^*, \xi), \\ 
\label{cts:3} 0 & = \beta \inn{\s_t}{\zeta} + (\nabla  \s, \nabla \zeta)  + \kappa(\sigma - \hat \sigma_B, \zeta)_{\Gamma} + (\lambda_c \hat h \s_1 + \lambda_c h_2 \s , \zeta) + B(\sigma - \hat \sigma_c, \zeta),  \\ 
\label{cts:4} 0 & = (\C (\E(\u) - \E^* \ph), \nabla \bet) - (\hat \g, \bet)_{\GN},
\end{alignat}
\end{subequations}
for a.e.~$t \in (0,T)$, and for all $\zeta, \xi \in H^1(\Omega)$ and $\bet \in \X$.

\subsection{Continuous dependence in weaker norms}
To recover the operator $\A$ in the first two equations,
we add to both sides of \eqref{cts:1} the term $(\mu, \zeta)$ and to both sides of \eqref{cts:2} the term $-\eps (\ph,\xi)$.  \andrea{Moreover, we define} a modified potential
\begin{align*}
\Psi(s) := \psi(s) - \frac{\eps}{2} s^2,
\end{align*} 
\andrea{which still} fulfils \eqref{Cts:pot:generalcase} as well as
\begin{align*}
\hat{\Psi}' := \Psi'(\ph_1) - \Psi'(\ph_2) = \hat{\psi}' - \eps \ph.
\end{align*}
In particular, \eqref{cts:1}, \andreaold{and} \eqref{cts:2} now \andrea{assume the form}
\begin{subequations}
\begin{alignat}{2}
\notag 0 & = \inn{\ph_t}{\zeta} + \inn{\A \mu}{ \zeta} - (\mu,\zeta) -\andrea{\Big(}\tfrac{\lambda_p}{1+|W_{,\E,1}|} ( \hat f \s_1 + f_2 \s)  - \lambda_a \hat k, \zeta\andrea{\Big)} \\
\label{cdn:1} & \quad - \Big (\tfrac{\lambda_p f_2 \s_2}{(1+|W_{,\E,1}|)(1+|W_{,\E,2}|)} (|W_{,\E,1}|-|W_{,\E,2}|), \zeta \Big ), \\
\label{cdn:2} 0 & = (\mu, \xi) - \eps \inn{\A \ph}{\xi} - \eps^{-1} (\hat{\Psi}', \xi) + \chi (\sigma, \xi) + (\C (\E(\u) - \E^* \ph): \E^*, \xi),
%\\
%\label{cdn:3} 0 & = \beta \inn{\s_t}{\zeta} + (\nabla \s,\andreaold{\nabla}  \zeta) + \kappa(\sigma - \hat \sigma_B, \zeta)_{\Gamma} + (\lambda_c \hat h \s_1 + \lambda_c h_2 \s , \zeta) + B(\sigma - \hat \sigma_c, \zeta)  , 
\end{alignat}
\end{subequations}
for a.e.~$t \in (0,T)$, and for all $\zeta, \xi \in H^1(\Omega)$.

First, testing \eqref{cts:3} with $\s$ gives
\begin{align}\label{cts:s}
\frac{\beta}{2} \frac{d}{dt} \no{\s}_{L^2}^2 + \no{\nabla \s}_{L^2}^2 + \frac{B}{2} \no{\s}_{L^2}^2 + \frac{\kappa}{2} \no{\s}_{\andrea{L^2_\Gamma}}^2 \leq C \Big ( \no{\ph}_{L^2}^2 + \no{\hat \s_c}_{L^2}^2 + \no{\hat \s_B}_{\andrea{L^2_\Gamma}}^2 \Big ),
\end{align}
where we neglected the non-negative term $\lambda_c h_2 \s^2$ and the constant $C$ is independent of $\beta$.  To this, we add the equalities obtained from testing \eqref{cdn:1} by $\A^{-1} \ph$, \eqref{cdn:2} with $\ph$ and \eqref{cts:4} with $\u$, whilst recalling the relation $\inn{\A \mu}{\A^{-1}\ph} = \inn{\mu}{\ph} = (\mu,\ph)$, and recalling the coercivity estimate from \eqref{prop:E:1}:
\begin{align}\label{C:coer}
\C \E(\u): \E(\u) \geq C_4 \no{\E(\u)}_{L^2}^2,
\end{align}
we arrive at
\begin{align*}
& \frac{1}{2} \frac{d}{dt} \Big(\no{\ph}_*+\beta\no{\s}_{L^2}^2 \Big)	+ \eps \no{\ph}_{H^1}^2 + \no{\nabla \s}_{\andreaold{L^2}}^2 + \frac{B}{2} \no{\s}_{L^2}^2 + \frac{\kappa}{2} \no{\s}_{\andrea{L^2_\Gamma}}^2 \\
& \qquad C_4 \no{\E(\u)}_{L^2}^2 + \eps^{-1} ( \hat{\Psi}', \ph)  - C \Big ( \no{\ph}_{L^2}^2 + \no{\hat \s_c}_{L^2}^2 + \no{\hat \s_B}_{\andrea{L^2_\Gamma}}^2 \Big ) \\
& \quad \leq (\C(\E(\u) - \E^* \ph):\E^* + \chi \s, \ph) + ( \mu,\andreaold{\A}^{-1} \ph) \\
& \qquad + \andrea{\Big(}\tfrac{\lambda_p}{1+|W_{,\E,1}|} (\hat f \s_1 + f_2 \s) - \lambda_a \hat k, \andreaold{\A}^{-1} \ph\andrea{\Big)} \\
& \qquad + \Big ( \tfrac{\lambda_p f_2 \s_2}{(1 + |W_{,\E,1}|)(1+|W_{,\E,2}|)}(|W_{,\E,1}| - |W_{,\E,2}|), \A^{-1} \ph \Big ) \\
& \qquad + (\C \E^* \ph, \E(\u)) + (\hat \g, \u)_{\GN} \\
& \quad =: I_1 + I_2 + I_3 + I_4.
\end{align*}
Invoking the boundedness and Lipschitz continuity of $f$ and $k$, the boundedness of $\s_1$ and $\s_2$, the trace theorem, Young's inequality, Korn's inequality and the interpolation inequality \eqref{interpol_one}, we can estimate the terms on the right-hand side as follows:
\begin{align*}
I_1 & \leq \frac{\delta_1}{4} \no{\E(\u)}_{L^2}^2 + \frac{\eps}{8} \no{\ph}_{H^1}^2 + \frac{\delta_2}{2} \no{\s}^2 
\andreaold{+C \norma{\ph}_*^2}
+ (\mu, \andreaold{\A}^{-1} \ph), \\
I_2 & \leq  C ( \no{\ph}_{L^2} + \no{\s}_{L^2} )\no{\A^{-1} \ph}_{L^2} \leq \frac{\delta_2}{2} \no{\s}_{L^2}^2 + \frac{\eps}{8} \no{\ph}_{H^1}^2 + C \no{\ph}_{*}^2, \\
I_3 & \leq C ( \no{\E(\u)}_{L^2} + \no{\ph}_{L^2}) \no{\A^{-1} \ph}_{L^2} \leq \frac{\delta_1}{4} \no{\E(\u)}_{L^2}^2 + \frac{\eps}{8} \no{\ph}_{H^1}^2 + C \no{\ph}_{*}^2, \\
I_4 & \leq \frac{\delta_1}{2} \no{\E(\u)}_{L^2}^2 + \frac{\eps}{8} \no{\ph}_{H^1}^2 + C \Big ( \no{\ph}_{*}^2 + \no{\hat \g}_{L^2(\GN)}^2 \Big )
\end{align*}
for $\delta_1, \delta_2 >0$ yet to be determined.  Choosing $\delta_1 = \frac{C_4}{2}$ then gives
\begin{equation}\label{cts:est:1}
\begin{aligned}
& \frac{d}{dt} \Big ( \no{\ph}_{*}^2 + \beta \no{\s}_{L^2}^2 \Big ) + \andreaold{\eps} \no{\ph}_{H^1}^2 + 2 \no{\nabla \s}_{L^2}^2  \\
& \qquad + B \no{\s}_{L^2}^2 + \kappa \no{\s}_{\andrea{L^2_\Gamma}}^2 + C_4 \no{\E(\u)}_{L^2}^2 + \andreaold{2}\eps^{-1} (\hat{\Psi}', \ph) \\
& \quad \leq C \Big ( \no{\ph}_{*}^2 + \no{\hat \s_c}_{L^2}^2 + \no{\hat \s_B}_{\andrea{L^2_\Gamma}}^2 + \no{\hat \g}_{L^2(\GN)}^2 \Big ) + 2 \delta_2 \no{\s}_{L^2}^2 + 2 (\mu, \A^{-1} \ph).
\end{aligned}
\end{equation}
Next, from testing \eqref{cts:2} with $\A^{-1} \ph$ we infer that 
\begin{align*}
2 |(\mu, \A^{-1} \ph)| \leq 2 \no{\ph}_* \Big (\eps \no{\ph}_{H^1} + \eps^{-1} \andreaold{\no{\hat{\Psi}'}_{*}}  + \chi \no{\s}_{L^2} + c_*e_* \no{\E(\u)}_{L^2}+ c_* e_*^2 \no{\ph}_{L^2} \Big ),
\end{align*}
where
\begin{align*}
c_* := \max_{1 \leq i,j,k,l \leq d} \abs{\C_{ijkl}}, \quad e_* := \max_{1 \leq i,j \leq d} \abs{\E^*_{ij}}. 
\end{align*}
By the continuous embedding $L^{6/5}(\Omega) \subset H^1(\Omega)'$ and
$H^1(\Omega)\subset L^6(\Omega)$, \eqref{Cts:pot:generalcase} with
exponent $q = 4$, and H\"older's inequality we obtain that
\begin{align*}
\andreaold{\no{\hat{\Psi}'}_{*}} \leq C \andreaold{\no{\hat{\Psi}'}}_{L^{6/5}} \leq C \Big ( 1 + \no{\ph_1}_{L^6}^{4} + \no{\ph_2}_{L^6}^4 \Big ) \no{\ph}_{L^6} \leq C \no{\ph}_{H^1}
\end{align*}
taking into account $\ph_1, \ph_2 \in L^\infty(0,T;H^1(\Omega))$.  Then, by Young's inequality
\begin{align}\label{cts:mu:est}
2 |(\mu, \A^{-1} \ph)| \leq \frac{\eps}{\andreaold{4}} \no{\ph}_{H^1}^2 + \delta_2 \no{\s}_{L^2}^2 + \frac{C_4}{2} \no{\E(\u)}_{L^2}^2 + C \no{\ph}_*^2.
\end{align}
Next, recalling the convex-concave decomposition of the potential $\psi$, we deduce that 
\begin{align*}
\andreaold{2}\eps^{-1} (\hat{\Psi}', \ph) \geq \andreaold{2}\eps^{-1} \andreaold{(\psi_2'(\ph_1) - \psi_2'(\ph_2), \ph)} - \andreaold{2}\eps \no{\ph}_{L^2}^2 \geq - C \andrea{\no{\ph}_{L^2}^2} \geq - \frac{\eps}{\andreaold{4}} \no{\ph}_{H^1}^2 - C \no{\ph}_*^2.
\end{align*}
Substituting the above and \eqref{cts:mu:est} into \eqref{cts:est:1} leads to the differential inequality
\begin{align*}
& \frac{d}{dt} \Big ( \no{\ph}_*^2 + \beta \no{\s}_{L^2}^2 \Big ) + \frac{\eps}{\andreaold{2}}\no{\ph}_{H^1}^2 + \frac{C_4}{2} \no{\E(\u)}_{L^2}^2 \\
& \qquad + 2 \no{\nabla \s}_{L^2}^2 + \frac{B}{2} \no{\s}_{L^2}^2 + \frac{\kappa}{2} \no{\s}_{\andrea{L^2_\Gamma}}^2 - 3 \delta_2 \no{\s}_{L^2}^2 \\
& \quad \leq C \Big ( \no{\ph}_*^2 + \no{\hat \s_c}_{L^2}^2 + \no{\hat \s_B}_{\andrea{L^2_\Gamma}}^2 + \no{\hat \g}_{L^2(\GN)}^2 \Big ),
\end{align*}
with a constant $C$ independent of $\beta$.  For the case $\beta > 0$,
we can move the last term on the left-hand side to the right hand side
and invoke Gronwall's inequality and
Korn's inequality to deduce \eqref{cts} aside for the estimate of
$\mu$.  For the case $\beta = 0$, at least one of $\{B, \kappa\}$ is
non-zero, and so choosing $\delta_2$ sufficiently small and possibly invoking
the generalised Poincar\'e inequality (if $B = 0$), we can absorb the
contribution $3\delta_2\|\sigma\|^2_{L^2}$ on the left hand side, we deduce via Gronwall's inequality also an
analogous estimate.  Then, to complete the proof, from \eqref{cts:2}
we infer
\begin{align*}
\int_0^T \no{\mu}_*^2 \dt \leq C \int_0^T \Big (\no{\ph}_{H^1}^2 + \no{\E(\u)}_{L^2}^2 + \no{\s}_{L^2}^2 \Big ) \dt,
\end{align*}
which yields the remaining $L^2(0,T;H^1(\Omega)')$ estimate for $\mu$.

\subsection{Continuous dependence in stronger norms}
\label{SEC_cdstrong}
Now, suppose the exponent $q$ in \eqref{Cts:pot:generalcase} is 2, testing \eqref{cts:1} with $\eps \ph$ and \eqref{cts:2} with $\mu$ yields upon summing
\begin{align*}
\frac{\eps}{2} \frac{d}{dt} \no{\ph}_{L^2}^2 + \no{\mu}_{L^2}^2 &  = \eps \andrea{\Big(}\tfrac{\lambda_p}{1+|W_{,\E,1}|} ( \hat f \s_1 + f_2 \s)  - \lambda_a \hat k, \ph\andrea{\Big)} \\
& \quad + \eps \Big (\tfrac{\lambda_p f_2 \s_2}{(1+|W_{,\E,1}|)(1+|W_{,\E,2}|)} (|W_{,\E,1}|-|W_{,\E,2}|), \ph \Big ) \\
& \quad + \eps^{-1} ( \hat \psi',\mu) \andrea{-} \chi (\s, \mu) \andrea{-} (\C (\E(\u) - \E^* \ph) : \E^*, \mu) \\
& =: J_1 + J_2 + J_3.
\end{align*}
Then, by Young's inequality, the Lipschitz continuity of $f$ and $k$,
as well as the boundedness of $\sigma_1$ and $\sigma_2$, we obtain that
\begin{align*}
J_1 & \leq C\no{\ph}_{L^2}^2 + C\no{\s}_{L^2}^2, \\
J_2 & \leq C \no{\E(\u)}_{L^2}^2 + C \no{\ph}_{L^2}^2, \\
J_3 & \leq C \Big ( 1 + \no{\ph_1}_{L^\infty}^4 + \no{\ph_2}_{L^\infty}^4 \Big ) \no{\ph}_{L^2}^2 + C \no{\s}_{L^2}^2 + C \no{\E(\u)}_{L^2}^2 + \frac{1}{2} \no{\mu}_{L^2}^2.
\end{align*}
Employing the Gagliardo--Nirenberg inequality \eqref{GN:3d}, as well as the regularities $\ph_1, \ph_2 \in L^\infty(0,T;H^1(\Omega)) \cap L^2(0,T;H^2(\Omega))$ we see that 
\begin{align*}
\eps \frac{d}{dt} \no{\ph}_{L^2}^2 + \no{\mu}_{L^2}^2 \leq C \Big ( 1 + \no{\ph_1}_{H^2}^2 + \no{\ph_2}_{H^2}^2 \Big ) \no{\ph}_{L^2}^2 + C \no{\E(\u)}_{L^2}^2 + C \no{\s}_{L^2}^2\andrea{.}
\end{align*}
Then, Gronwall's inequality as well as the estimate \eqref{cts} for $\u$ and $\s$ yields the continuous dependence results for $\ph$ in $L^\infty(0,T;L^2(\Omega))$ and for $\mu$ in $L^2(0,T;L^2(\Omega))$.  Next, testing \eqref{cts:4} with $\andrea{\bet} = \u$ and using the coercivity estimate \eqref{C:coer}, as well as the continuous dependence result for $\ph$ in $L^\infty(0,T;L^2(\Omega))$, allows us to infer that 
\begin{align*}
\no{\E(\u)}_{L^2}^2 \leq C \no{\ph}_{L^2}^2 + C\no{\hat \g}_{L^2(\GN)}^2.
\end{align*}
Taking supremum in time and applying Korn's inequality yields the continuous dependence result for $\u$ in $L^\infty(0,T;\X)$.

\andreaold{
\subsection{Continuous dependence via time discretisation approach}
In this section, we assume \eqref{Cts:pot:generalcase} with exponent
\andrea{$q$ equal to $2$} and \eqref{Cts:last} to hold.  Recalling the
notation for the difference between two sets of weak solutions
\eqref{diff_one}-\eqref{diff_two}, for $t \in (0,T)$ and $h > 0$, we
define
\begin{align*}
	\d_h \ph(t):= \frac {\ph(t)-\ph(t-h)} h, 
	\quad
	\d_h \s(t):= \frac {\s(t)-\s(t-h)}{h},
	\quad 
	\d_h \u(t):= \frac {\u(t)-\u(t-h)} h.
\end{align*}
Moreover, we set
\begin{align*}
	\ph(t):= \ph_{0,1} - \ph_{0,2}, 
	\quad
	\s(t):= \s_{0,1} - \s_{0,2}, 
	\quad
	\u(t):= \u_{0,1} - \u_{0,2},
	\quad \text{ for all } t \leq 0,
\end{align*}
where $\u_{0,i}$, $i \in \{1,2\}$, is the unique solution of the following elliptic problem:
\begin{align*}
	(W_{,\E}(\ph_{0,i}, \E(\u_{0,i})), \nabla \bet) = (\C(\E(\u_{0,i}) - \andrea{\hat \E - \E^* \ph_{0,i}}), \nabla \bet) = (\g_i, \bet)_{\GN} 
		\quad \forall\, \bet \in \X.
\end{align*}
The coercivity estimate \eqref{C:coer} and the Lax--Milgram theorem
gives the existence of a unique $\u_{0,i} \in \X$, and by taking the
difference of the elliptic equations for $\u_{0,1}$ and $\u_{0,2}$ we
can infer from testing $\andrea{\bet} = (\u_{0,1} - \u_{0,2})$ and Korn's inequality that
\begin{align}\label{u0:est}
\no{\u_{0,1} - \u_{0,2}}_{H^1}^2 \leq C \no{\ph_{0,1} - \ph_{0,2}}_{L^2}^2 + C \no{\hat \g}_{L^2(\GN)}^2.
\end{align}
Let us point out that the difference $\hat U := U(\ph_1,\s_1, \E(\u_1))-U(\ph_2,\s_2, \E(\u_2))$ is already bounded in $L^2(Q)$ by previous results. Namely from \eqref{cdn:1}, we have
\begin{align*}
 \no{\hat U}_{L^2(Q)} &
	\leq C (\no{\ph}_{L^2(Q)} + \no{\s}_{L^2(Q)} +\no{\E(\u)}_{L^2(Q)} )
	\\ 
	&  \leq 
	C\big (\| \ph_{0,1} - \ph_{0,2} \|_{L^2}
	+\beta \| \s_{0,1} - \s_{0,2} \|_{L^2} 
	+ \no{\hat \g}_{L^2(\GN)} + \no{\hat \s_{B}}_{L^2(\Sigma)}+ \no{\hat \s_{c}}_{L^2(Q)} \big ) \\
	& =: C \mathcal{Y},
\end{align*}
with a constant $C$ independent of $\beta$.  Similarly, for $\hat S := S(\ph_1,\s_1)- S(\ph_2,\s_2)$, we have
\begin{align*}
 \no{ \hat S}_{L^2(Q)}  & \leq C (\norma{\ph}_{L^2(Q)} + \no{\s}_{L^2(Q)}  ) \leq C\mathcal{Y},
%	\\ & \leq 
%	C\big (\| \ph_{0,1} - \ph_{0,2} \|_{L^2}
%	+\beta \| \s_{0,1} - \s_{0,2} \|_{L^2} 
%	+ \no{\hat \g}_{L^2(\GN)} + \no{\hat \s_{B}}_{L^2(\Sigma)}+ \no{\hat \s_{c}}_{L^2(Q)} \big ), 
\end{align*}
with a constant $C$ independent of $\beta$.
%where
%\begin{align*}
%	S:= S(\ph_1,\s_1)- S(\ph_2,\s_2) 
%	= - \lambda_c \hat h \s_1 - \lambda_c h_2 \s
%		+ B (\s-\hat \s_c).
%\end{align*}
%Moreover, let us underline that we are now working with the difference of solutions
%so that by the above extension we realise that
%\begin{align*}
%	\ph(t)=0, \quad 
%	\s(t)=0, \quad 
%	\u(t)=0, \quad \hbox{for all $t\leq 0$}.
%\end{align*}
}

Without loss of generality, we fix $t\in (0,T)$, and take $h$ sufficiently small so that $t- h > 0$.  Then, integrating \eqref{cts:1} over time from $t-h$ to $t$, and dividing by $h$, we have
\begin{align}\label{cts:1:time}
0 = (\d_h \ph(t),\zeta) + \frac{1}{h}
\andrea{\Big(} \int_{t-h}^t \nabla \mu(\tau) d \tau, \nabla \zeta \andrea{\Big)}
- \frac{1}{h} \andrea{\Big(} \int_{t-h}^t \hat U(\tau) d \tau, \zeta \andrea{\Big)}
\end{align}
holding for all $\zeta \in H^1(\Omega)$.  Choosing $\zeta = \mu(t)$,
and combining the resulting equality with the one obtained from
testing \eqref{cts:2} with $\xi = \d_h \ph(t)$, we note \andrea{that}
a cancellation occurs and obtain
\begin{align*}
0 & = \eps (\nabla \ph(t), \nabla \d_h \ph(t)) + \eps^{-1}(\hat \psi'(t), \d_h \ph(t)) - \chi (\s(t), \d_h \ph(t))\\
& \quad - (\C(\E(\u(t)) - \E^* \ph(t)): \E^*, \d_h \ph(t)) + \frac{1}{h} \int_{t-h}^t  (\nabla \mu(\tau), \nabla \mu(t)) - (\hat U(\tau), \mu(t)) d \tau.
\end{align*}
To the above, we add the equality obtained from testing
$\bet = \d_h \u (t)$ in \eqref{cts:4}, as well as using the identity from
testing $\zeta = \eps^{-1} \hat \psi'(t) -\chi \s(t)$ in
\eqref{cts:1:time}, leading to
\begin{equation}\label{ctsdep:time:1}
\begin{aligned}
0 & = \eps (\nabla \ph(t), \nabla \d_h \ph(t)) + (\C(\E(\u(t)) - \E^* \ph(t)), \d_h (\E(\u(t))- \E^*\ph(t))) \\
& \quad  - (\hat \g, \d_h \u(t))_{\GN} + \frac{1}{h} \int_{t-h}^t (\nabla \mu(\tau), \nabla (\mu(t) - \eps^{-1} \hat \psi'(t) + \chi \s(t)) d\tau \\
& \quad  - \frac{1}{h} \int_{t-h}^t (\hat U(\andrea{\tau}), \mu(t) - \eps^{-1} \hat \psi'(t) + \chi \s(t)) d\tau.
\end{aligned}
\end{equation}
We aim to send $h \to 0$ rigorously.  By the identity $a(a-b) = \frac{1}{2}(a^2 - b^2 + (a-b)^2)$ and the Lebesgue differentiation theorem, we see that for any $s > 0$,
\begin{align*}
& \int_0^s \eps (\nabla \ph(t), \nabla \d_h \ph(t)) \dt \\
& \quad = \frac{\eps}{2h} \int_0^s \no{\nabla \ph(t)}_{L^2}^2 - \no{\nabla \ph(t-h)}_{L^2}^2 + \no{\nabla (\ph(t) - \ph(t-h))}_{L^2}^2 \dt \\
& \quad \geq \frac{\eps}{2h} \Big ( \int_{s-h}^s \no{\nabla \ph(t)}_{L^2}^2 \dt \andrea{-} \int_{-h}^0 \no{\nabla \ph(t)}_{L^2}^2 \dt \Big ) \to \frac{\eps}{2} \no{\nabla \ph(s)}_{L^2}^2 - \frac{\eps}{2} \no{\nabla \ph(0)}_{L^2}^2.
\end{align*}
Similarly, a short calculation using the symmetry of the constant elasticity tensor $\C$ shows that
\begin{align*}
& \int_0^s (\C(\E(\u(t)) - \E^* \ph(t)), \d_h (\E(\u(t)) - \E^* \ph(t))) \dt \\
& \quad =\frac{1}{h} \int_0^s \int_\Omega W(\ph(t), \E(\u(t))) - W(\ph(t-h), \E(\u(t-h))) d\x \dt \\
& \qquad + \frac{1}{2h} \int_0^s \no{\sqrt{\C} ( \E(\u(t)- \u(t-h)) - \E^* (\ph(t) - \ph(t-h)))}_{L^2}^2 \dt \\
& \quad \geq \frac{1}{h} \int_{s-h}^s \int_\Omega W(\ph(t),\E(\u(t))) d\x \dt - \frac{1}{h} \int_{-h}^0 \int_\Omega W(\ph(t), \E(\u(t))) d\x \dt \\
& \quad \to \int_\Omega W(\ph(s), \E(\u(s))) d\x - \int_\Omega W(\ph(0), \E(\u(0))) d\x.
\end{align*}
Meanwhile, 
\begin{align*}
\int_0^s ( \hat \g, \d_h \u(t))_{\GN} \dt & = \frac{1}{h} \int_{s-h}^s ( \hat \g, \u(t))_{\GN} \dt - \frac{1}{h} \int_{-h}^0 (\hat \g, \u(t))_{\GN} \dt \\
& \to (\hat \g, \u(s))_{\GN} - (\hat \g, \u(0))_{\GN},
\end{align*}
and
\andrea{\begin{align*}
\frac{1}{h} \int_{s-h}^s (\nabla \mu(t), \nabla \zeta) - (\hat U(t), \zeta) \dt \to (\nabla \mu(s), \nabla \zeta) - (\hat U(s), \zeta) \quad \forall\, \zeta \in H^1(\Omega).
\end{align*}}
Hence, integrating \eqref{ctsdep:time:1} over $(0,s)$ and passing to the limit $h \to 0$ yields 
\begin{align*}
& \frac{\eps}{2} \no{\nabla \ph(s)}_{L^2}^2 + \int_\Omega W(\ph(s), \E(\u(s))) d\x + \int_0^s \no{\nabla \mu(t)}_{L^2}^2 \dt \\
& \quad \leq \frac{\eps}{2} \no{\nabla \ph(0)}_{L^2}^2 + \no{W(\ph(0), \E(\u(0)))}_{L^1} + (\hat \g, \u(0) - \u(s))_{\GN} \\
& \qquad + \int_0^s (\nabla \mu(t), \eps^{-1} \nabla \hat \psi'(t) - \chi \nabla \s(t)) + (\hat U(t), \mu(t) - \eps^{-1} \hat \psi'(t) + \chi \s(t)) \dt.
\end{align*}
Then, by Young's inequality, \eqref{cts} and \eqref{u0:est}, the right-hand side can be bounded by
\begin{equation}\label{ctsdep:time:2}
\begin{aligned}
\mathrm{RHS} & \leq C \Big (  \no{\ph(0)}_{H^1}^2 + \no{\u(0)}_{H^1}^2 + \no{\hat \g}_{L^2(\GN)}^2 + \no{\u(s)}_{H^1}^2 \Big ) + \frac{1}{2} \int_0^s \no{\nabla \mu(t)}_{L^2}^2 \dt  \\
& \quad + C \no{\s}_{L^2(0,T;H^1)}^2 + C \no{\hat U}_{L^2(Q)}^2 + C \no{\mu}_{L^2(Q)}^2 + C\no{\hat \psi'}_{L^2(0,T;H^1)}^2 \\
 & \leq C \mathcal{Y}^2 +  \frac{1}{2} \int_0^s \no{\nabla \mu(t)}_{L^2}^2 \dt + C\no{\hat \psi'}_{L^2(0,T;H^1)}^2.
\end{aligned}
\end{equation}
As a consequence of the calculations in Section \ref{SEC_cdstrong} and the continuous dependence estimate \eqref{cts}, we have
\begin{align}\label{ctsdep:time:3}
\no{\hat \psi'}_{L^2(Q)}^2 \leq C \Big ( 1 +\no{\ph_1}_{L^2(0,T;H^2)}^2 + \no{\ph_2}_{L^2(0,T;H^2)}^2 \Big ) \no{\ph}_{L^\infty(0,T;L^2)}^2 \leq C \mathcal{Y}^2,
\end{align}
while invoking \eqref{Cts:pot:generalcase} and \eqref{Cts:last} we see that 
\begin{align*}
\no{\nabla \hat \psi'}_{L^2(Q)}^2 & \leq C \no{(\psi''(\ph_1) - \psi''(\ph_2)) \nabla \ph_1}_{L^2(Q)}^2 + C \no{\psi''(\ph_2) \nabla \ph}_{L^2(Q)}^2 \\
& \leq C \int_Q (1 + |\ph_1|^2 + |\ph_2|^2)|\ph|^2 |\nabla \ph_1|^2  + (1 + |\ph_2|^4) |\nabla \ph|^2 d\x \dt \\
& \leq C \int_0^T ( 1 + \no{\ph_1}_{L^6}^2 + \no{\ph_2}_{L^6}^2)\no{\ph}_{L^6}^2 \no{\nabla \ph_1}_{L^6}^2  + (1 + \no{\ph_2}_{H^2}^2) \no{\nabla \ph}_{L^2}^2 \dt \\
& \leq C \int_0^T ( 1 + \no{\ph_1}_{H^2}^2 + \no{\ph_2}_{H^2}^2 ) \no{\ph}_{H^1}^2 \\
& \leq C \mathcal{Y}^2 + C\int_0^T ( 1 + \no{\ph_1}_{H^2}^2 + \no{\ph_2}_{H^2}^2 ) \no{\nabla \ph}_{L^2}^2,
\end{align*}
where we have used the second inequality of \eqref{ctsdep:time:3}.  Hence, from \eqref{ctsdep:time:2}, we infer that
\begin{align*}
\no{\nabla \ph(s)}_{L^2}^2 + \int_0^s \no{\nabla \mu(t)}_{L^2}^2 \dt \leq C\mathcal{Y}^2 +  C\int_0^T ( 1 + \no{\ph_1}_{H^2}^2 + \no{\ph_2}_{H^2}^2 ) \no{\nabla \ph}_{L^2}^2,
\end{align*}
and the result follows first from the application of the integral form of Gronwall's inequality and then the observation that by elliptic regularity
\begin{align*}
\no{\ph}_{L^2(0,T;H^2)}^2 & \leq C \Big ( \no{\ph}_{L^2(0,T;H^1)}^2 + \no{\mu}_{L^2(Q)}^2 + \no{\hat \psi'}_{L^2(Q)}^2 \Big ) \\
& \quad + C \Big ( \no{\s}_{L^2(Q)}^2 + \no{\C(\E(\u) - \E^* \ph): \E^*}_{L^2(Q)}^2 \Big ) \\
& \leq C \mathcal{Y}^2.
\end{align*}
For the nutrient, under the hypothesis on the data
$\sigma_{B_1},\sigma_{B_2} \in H^1(0,T;L^2(\Gamma))$, we infer from
Theorem \ref{thm:Reg} that
\andrea{$\sigma_1,\sigma_2 \in H^1(0,T;L^2(\Omega))) \cap
  L^\infty(0,T;H^1(\Omega))$}.  Then, in \eqref{cts:3}
the term $\beta \inn{\s_t}{\zeta}$ can be written as
$\beta (\s_t, \zeta)$ and choosing $\zeta= \d_h \s(t)$ and
integrating over $(0,s)$ yields with the help of previous calculations
\begin{align*}
0 & =\int_0^s  (\beta \s_t(t) + \hat S(t)
\andrea{+ B(\s(t)-\hat \s_c(t))}, \d_h \s(t)) + (\nabla \s(t), \nabla \d_h \s(t)) 
\\ & \quad + \kappa (\s(t) - \hat \s_B(t), \d_h \s(t))_{\Gamma} \dt \\
& \geq \andrea{\frac{1}{2h}} \int_0^s  \Big (\no{\nabla
  \s(t)}_{L^2}^2+B\|\sigma(t)\|^2_{L^2} + \kappa
  \no{\s(t)}_{\andrea{L^2_\Gamma}}^2 
  \\ & \hspace{2cm} 
  - \no{\nabla \s(t-h)}_{L^2}^2 -B\|\sigma(t-h)\|^2_{L^2} - \kappa \no{\s(t-h)}_{\andrea{L^2_\Gamma}}^2 \Big ) \dt \\
& \quad + \int_0^s  (\beta\s_t(t) + \hat S(t)-B\hat{\sigma}_c(t), \d_h
  \s(t)) - \kappa (\hat \s_B(t), \d_h \s(t))_{\Gamma} \dt.
\end{align*}
By a change of variables we deduce that 
\begin{align*}
0 & \geq \andrea{\frac{1}{2h}} \int_{s-h}^s \Big (\no{\nabla \s(t)}_{L^2}^2+B\|\sigma(t)\|^2_{L^2} + \kappa \no{\s(t)}_{\andrea{L^2_\Gamma}}^2 \Big ) \dt 
\\& \quad 
- \andrea{\frac{1}{2h}}  \int_{-h}^0 \Big ( \no{\nabla \s(t)}_{L^2}^2+B\|\sigma(t)\|^2_{L^2} + \kappa \no{\s(t)}_{\andrea{L^2_\Gamma}}^2 \Big ) \dt \\
& \quad + \int_0^s  (\beta\s_t(t) + \hat S(t) \andrea{- B \hat \s_c(t)}, \d_h \s(t)) \dt + \int_0^s (\d_h \hat \s_B(t), \s(t))_{\Gamma} \dt \\
& \quad - \frac{1}{h} \int^s_{s-h}\big( \hat \s_B(t), \s(t)\andrea{{\big)}_\Gamma} \dt + \frac{1}{h}\int_{-h}^0 \big(\hat \s_B(t+h), \s(t)\andrea{\big)_{\Gamma}\dt}.
\end{align*}
Setting $\hat \s_B(t) = \hat \s_B(0)$ for $t \leq 0$, and applying the Lebesgue differentiation theorem, we find that in the limit $h \to 0$ it holds
\begin{align*}
& \frac{1}{2} \Big ( \no{\nabla \s(s)}_{L^2}^2 \andrea{+ B \norma{\s(s)}^2_{L^2}} + \kappa \no{\s(s)}_{\andrea{L^2_\Gamma}}^2 \Big ) + \andrea{\frac{\beta}{2}} \int_0^s \no{\s_t}_{L^2}^2 \dt \\
& \quad \leq \frac{1}{2} \Big ( \no{\nabla \s(0)}_{L^2}^2 \andrea{+B \norma{\s(0)}^2_{L^2}}+ \kappa \no{\s(0)}_{\andrea{L^2_\Gamma}}^2 \Big ) \andrea{+C\norma{\hat\s_c}^2_{L^2(Q)}} + \andrea{C} \no{\hat S}_{L^2(Q)}^2  \\
& \qquad + C \no{\hat \s_{B,t}}_{L^2(\Sigma)}^2 + C \no{\hat \s_{B}}_{C^0([0,T]; L^2(\Gamma))}^2 + C \no{\s}_{\andrea{C^0([0,T];L^2(\Gamma))}}^2  \\
& \quad \leq C \mathcal{Y}^2 + C \no{\s(0)}_{H^1}^2 + C \no{\hat \s_B}_{H^1(0,T;L^2(\Gamma))}^2.
\end{align*}
Furthermore, \andrea{assuming $\s_{B_1},\s_{B_2} \in L^2(0,T;H^{1/2}(\Gamma))$,} elliptic regularity gives
\begin{align*}
& \no{\s}_{L^2(0,T;H^2)}^2 
\\
 & \quad \leq C \Big ( \no{\s_t}_{L^2(Q)}^2 + \no{\s}_{L^2(0,T;H^1)}^2 \andrea{+\norma{\hat\s_c}^2_{L^2(Q)}} + \no{\hat S}_{L^2(Q)}^2 + \no{\hat \s_B}_{L^2(0,T;H^{1/2}(\Gamma))}^2 \Big ),
\end{align*}
and this completes the proof.

\section*{Acknowledgement}
The third author gratefully acknowledges 
financial support from the LIA-COPDESC initiative and from the research training group  2339 ``Interfaces,
Complex Structures, and Singular Limits'' of the German Science Foundation (DFG).

\end{document}